\documentclass[final,5p,times,twocolumn]{elsarticle}

\usepackage{amssymb}
\usepackage{amsmath}
\usepackage{bbm}
\usepackage{graphicx}
\usepackage{exscale}
\usepackage{clrscode}
\usepackage{hhline}

\biboptions{sort&compress}

\newcommand{\E}[1]{\mathop{{\rm \bf E}\!\left\{#1\right\}}\nolimits}

\newdefinition{dfn}{Definition}

\newtheorem{exampe}{Example}[section]

\begin{document}

\begin{frontmatter}
\title{MATLAB-based general approach for square-root extended-unscented and fifth-degree cubature Kalman filtering methods}

\author[CEMAT]{M.~V.~Kulikova\corref{cor}} \ead{maria.kulikova@ist.utl.pt} \cortext[cor]{Corresponding
author.}

\author[CEMAT]{G.~Yu.~Kulikov} \ead{gkulikov@math.ist.utl.pt}

\address[CEMAT]{CEMAT, Instituto Superior T\'ecnico, Universidade de Lisboa, Av.~Rovisco Pais, 1049-001 Lisboa, Portugal.}

\begin{abstract}
A stable square-root approach has been recently proposed for the unscented Kalman filter (UKF) and fifth-degree cubature Kalman filter (5D-CKF) as well as for the mixed-type methods consisting of the extended Kalman filter (EKF) time update and the UKF/5D-CKF measurement update steps. The mixed-type estimators provide a good balance in trading between estimation accuracy and computational demand because of the EKF moment differential equations involved. The key benefit is a consolidation of reliable state mean and error covariance propagation by using delicate discretization error control while solving the EKF moment differential equations and an accurate measurement update according to the advanced UKF and/or 5D-CKF filtering strategies. Meanwhile the drawback of the previously proposed estimators is an utilization of sophisticated numerical integration scheme with the built-in discretization error control that is, in fact, a complicated and computationally costly tool. In contrast, we design here the mixed-type methods that keep the same estimation quality but reduce a computational time significantly.  The novel estimators elegantly utilize any MATLAB-based numerical integration scheme developed for solving ordinary differential equations (ODEs) with the required accuracy tolerance pre-defined by users. In summary, a simplicity of the suggested estimators, their numerical robustness with respect to roundoff due to the square-root form utilized as well as their estimation accuracy due to the MATLAB ODEs solvers with discretization error control involved are the attractive features of the novel estimators.  The numerical experiments are provided for illustrating a performance of the suggested methods in comparison with the existing ones.
\end{abstract}


\end{frontmatter}

\section{Introduction}

In this paper, we focus on Bayesian filtering methods for estimating the hidden state of nonlinear continuous-discrete stochastic systems that imply  a finite number of summary statistics to be calculated instead of dealing with entire posterior density function information. The most famous techniques of such sort are based on computation of the first two moments, i.e. the mean and covariance matrix, and include the extended Kalman filter (EKF) strategy~\cite{1964:Ho,1970:Jazwinski:book}, the Unscented Kalman filtering (UKF) approach~\cite{2000:Julier,2001:Wan:book,2004:Julier,2015:Menegaz}, the quadrature Kalman filtering (QKF) methods~\cite{2007:Arasaratnam,2000:Ito}, the third-degree Cubature Kalman filtering (3D-CKF) algorithms~\cite{2009:Haykin,2010:Haykin} as well as the recent fifth-degree Cubature Kalman filtering (5D-CKF) suggested in~\cite{2013:Automatica:Jia,2018:Haykin}. As indicated in~\cite{2009:Haykin}, these sub-optimal estimators suggest a good practical feasibility, attractive simplicity in their implementation and good estimation quality with the fast execution time. All of these makes such estimators the methods of choice when solving practical applications.

The {\it continuous-discrete} nonlinear Kalman-like filtering methods might be designed under two alternative frameworks~\cite{2012:Frogerais,2014:Kulikov:IEEE}. The first approach implies the use of numerical schemes for solving the given stochastic differential equation (SDEs) of the system at hand. It is usually done by using either Euler-Maruyama method or It\^{o}-Taylor expansion for discretizing the underlying SDEs; see~\cite{1999:Kloeden:book}. More precisely, the readers may find such continuous-discrete filtering algorithms designed for the EKF in~\cite{1970:Jazwinski:book}, the 3D-CKF in~\cite{2010:Haykin,2020:Automatica:Kulikova}, the 5D-CKF in~\cite{2018:Haykin}, the UKF in~\cite{2019:Leth} and in many other studies. However the main drawback of the mentioned implementation framework is that SDEs solvers utilize prefixed and usually equidistant discretization meshes set by users, i.e. the resulting filtering methods require a special manual tuning from users prior to filtering in order to find an appropriate mesh that ensures a small discretization error arisen while propagating the mean and filter covariance matrix. Clearly, it is hardly predictable situation and this makes the resulting continuous-discrete estimators to be incapable to manage the missing measurement cases and they are inaccurate when solving estimation problems with irregular and/or long sampling intervals. Finally, any SDEs solver implies no opportunity to control and/or bound the occurred discretization error due to its stochastic nature. An alternative implementation framework assumes the derivation of the related filters' moment differential equations and then an utilization of the numerical methods derived for solving ordinary differential equations (ODEs); e.g., see the discussion in~\cite{2012:Frogerais,2014:Kulikov:IEEE,2018:Haykin}. This approach leads to a more accurate implementation framework because the discretization error control techniques are possible to involve and this makes the resulting filtering algorithms to be self-adaptive. More precisely, to perform the propagation step in the most accurate way, we control and keep the discretization error less than the pre-defined tolerance value given by users. The latter implies that the filtering algorithms of such kind enjoy self-adaptation mechanisms, that is, the discretization meshes in use are generated automatically by the filters themselves and with no user's effort. The user needs only to restrict the maximum magnitude of the discretization error tolerated in the corresponding state estimation run. It is worth noting here that the moment differential equations are derived for the continuous-discrete UKF and 3D-CKF estimators in~\cite{2007:Sarkka,2012:Sarkka} as well as for the EKF in~\cite{1970:Jazwinski:book}. To the best of the authors' knowledge, the 5D-CKF moment differential equations are not derived, yet. It could be an area for a future research.

In this paper, we focus on the mixed-type estimators consisting of the EKF time update and the UKF/5D-CKF measurement update steps. The key benefit of such methods is an accurate state mean and error covariance propagation combined with accurate measurement updates with a good balance in trading between estimation accuracy and computational demand. The boosted estimation quality at the filters' propagation step owes it all to the implemented discretization error control when solving the EKF moment differential equations involved.  More precisely, we have previously proposed such advanced mixed-type EKF-UKF filter and its stable square-root implementations in~\cite{2017:SP:Kulikov,2020:SP:Kulikov} as well as the EKF-5DCKF method with its square-rooting solution in~\cite{2020:IETSonar:KulikovaKulikov}. In the previously published works, it is performed by using the variable-step size nested implicit Runge-Kutta formulas with the built-in automatic global error control, which is a quite complicated and computationally costly tool for practical use~\cite{2013:Kulikov:IMA}. In contrast to the existing results, in this paper we propose the mixed-type estimation methods with the same improved estimation quality but with the reduced computational time and simple implementation way. The novel estimators elegantly utilize any MATLAB-based numerical integration scheme developed for solving ODEs with the required accuracy tolerance pre-defined by users. Additionally, we suggest a general square-rooting implementation scheme that is valid for both the UKF and CKF strategies and provides inherently more stable algorithms than the conventional filtering implementations; see more details about square-root methods in~\cite[Chapter~5]{1977:Bierman:book}, \cite[Chapter~12]{2000:Kailath:book}, \cite[Section~5.5]{2006:Simon:book}, \cite[Chapter~7]{2015:Grewal:book} and many other studies. Further advantages of our novel filtering methods are (i) an accurate processing of the missing measurements and/or irregular sampling interval scenario in an automatical mode, and (ii) the derivative free UKF and 5D-CKF strategies involved are of particular interest for estimating the nonlinear stochastic systems with non-differentiable observation equations~\cite{2017:Kulikov:IEEE:mix}. In summary, a simplicity of the suggested estimators, their numerical robustness in a finite precision arithmetics due to the square-root form as well as their estimation accuracy due to the discretization error control provided by MATLAB ODEs solver involved are the attractive features of the novel mixed-type filters. All this makes them to be preferable filtering implementation methods for solving practical applications. The numerical experiments are also provided for illustrating a performance of the suggested methods compared to the counterparts existing in engineering literature.

The paper is organized as follows. Section~\ref{problem:statement} provides an overview of the continuous-discrete EKF-UKF and EKF-5DCKF methods under examination. Section~\ref{main:result} suggests the new result that is a general MATLAB-based square-root scheme applied to the mentioned estimators. Finally, the results of numerical experiments are discussed in Section~\ref{numerical:experiments} and then Section~\ref{conclusion} concludes the paper by a brief overview of the problems that are still open in this research domain.

\section{Continuous-discrete EKF-UKF and EKF-5DCKF} \label{problem:statement}

Consider continuous-discrete stochastic system of the form
\begin{align}
dx(t) & = f\bigl(t,x(t)\bigr)dt+Gd\beta(t), \quad t>0,  \label{eq1.1} \\
z_k   & =  h(k,x(t_{k}))+v_k, \quad k =1,2,\ldots \label{eq1.2}
\end{align}
%The process model in the stochastic system (\ref{eq1.1}) and (\ref{eq1.2}) is an It$\hat{\rm o}$-type {\em Stochastic Differential Equation} (SDE).
where  $x(t) \in \mathbb R^{n}$ is the unknown state vector to be estimated and  $f:\mathbb R\times\mathbb
R^{n}\to\mathbb R^{n} $ is the time-variant drift function. The process uncertainty is modelled by the additive noise term where
$G \in \mathbb R^{n\times q}$ is the time-invariant diffusion matrix and $\beta(t)$ is the $q$-dimensional Brownian motion whose increment $d\beta(t)$ is independent of $x(t)$ and has the covariance $Q\,dt>0$. Finally, the available data $z_k = z(t_{k})$ is measured at some discrete-time points $t_k$, and $z_k \in \mathbb R^{m}$ comes with the sampling rate (sampling period) $\Delta_k=t_{k}-t_{k-1}$. The measurement noise term $v_k$ is assumed to be white Gaussian noise with zero mean and known covariance $R_k>0$, i.e. $\E{v_kv_j^{\top}} = R_k\delta_{kj}$ and $\delta_{kj}$ is the Kronecker delta function. Finally, the initial state $x(t_0) \sim {\cal N}(\bar x_0,\Pi_0)$ with $\Pi_0>0$, and the noise processes are all assumed to be statistically independent.

The first part of the novel mixed-type filtering methods is the {\it continuous-discrete} EKF implementation framework. Following~\cite{1970:Jazwinski:book}, the Moment Differential Equations (MDEs) should be solved in each sampling interval $[t_{k-1},t_k]$ with the initial $\hat x(t_{k-1}) =\hat x_{k-1|k-1}$ and $P(t_{k-1})=P_{k-1|k-1}$ as follows:
\begin{align}
    \frac{d\hat x(t)}{dt} & = f\bigl(t,\hat x(t)\bigr),  \label{eq2.1} \\
    \frac{dP(t)}{dt}& = F\bigl(t,\hat x(t)\bigr) P(t)+P(t)F^{\top}\bigl(t,\hat x(t)\bigr)+ GQG^{\top}  \label{eq2.2}
 \end{align}
where $F\bigl(t,\hat x(t)\bigr)={\partial_{x} f\bigl(t,x(t)\bigr)}|_{x(t)=\hat x(t)}$ is the Jacobian matrix.

In this paper, we suggest a general implementation framework for solving~\eqref{eq2.1}, \eqref{eq2.2} in an accurate way through the existing MATLAB's  ODEs solvers with the involved discretization error control; see the full list in Table~12.1 in~\cite[p.~190]{2005:Higham:book}. We stress that any method can be utilized in our novel mixed-type algorithms and this general MATLAB-based implementation scheme is proposed in Appendix~A.

The mixed-type filter consolidates the suggested accurate MATLAB-based EKF time update step with the advanced 5D-CKF rule for processing the new measurement and for computing the filtered ({\it a posteriori}) estimate $\hat x_{k|k}$. Following~\cite{2018:Haykin}, the fifth-degree spherical-radial cubature rule is utilized for approximating the Gaussian-weighted integrals by the formula
\begin{align*}
&\int \limits_{\mathbb R^n} f(x){\mathcal N}(x;\mu,\Sigma) dx \approx  \frac{2}{n+2}f(\mu)+\frac{1}{(n+2)^2} \\
\times&\biggl[\sum \limits_{j=1}^{n(n-1)/2} \left\{ f(\mu+\sqrt{\Sigma}\sqrt{n+2}s_j^+) + f(\mu-\sqrt{\Sigma}\sqrt{n+2}s_j^+)\right\}  \\
+&\sum \limits_{j=1}^{n(n-1)/2} \left\{ f(\mu+\sqrt{\Sigma}\sqrt{n+2}s_j^-) + f(\mu-\sqrt{\Sigma}\sqrt{n+2}s_j^-)\right\} \biggr] \\
+&\frac{4-n}{2(n+2)^2}\!\sum \limits_{i=1}^{n}\left\{ f(\mu+\sqrt{\Sigma}\sqrt{n+2}e_i) + f(\mu-\sqrt{\Sigma}\sqrt{n+2}e_i)\right\}.
\end{align*}
The matrix $\Sigma = \sqrt{\Sigma}\sqrt{\Sigma}^{\top}$ and $\sqrt{\Sigma}$ is its square-root factor where $e_i$ denotes the $i$-th unit coordinate vector in $\mathbb R^n$ with
\begin{align}
\{s^+_j\} & = \left\{ \sqrt{\frac{1}{2}}(e_k+e_i): k<l, k,l = 1,2,\ldots, n\right\}, \label{eq:splus} \\
\{s^-_j\} & = \left\{ \sqrt{\frac{1}{2}}(e_k-e_i): k<l, k,l = 1,2,\ldots, n\right\}. \label{eq:sminus}
\end{align}

Denote the 5D-CKF points and weights utilized as follows:
\begin{align}
w_0     & = \frac{2}{n+2}           &  \xi_0    & = \mathbf{0}_n, \label{5dckf:w0} \\
w_j^+   & = \frac{1}{2(n+2)^2},     &  \xi^+_j  & = \sqrt{n+2}s^+_j, & \xi^+_{l+j} & = -\sqrt{n+2}s^+_j, \label{5dckf:wj1} \\
w_j^-   & = \frac{1}{2(n+2)^2},     &  \xi^-_j  & = \sqrt{n+2}s^-_j, & \xi^-_{l+j} & = -\sqrt{n+2}s^-_j, \label{5dckf:wj2}\\
w_p     & = \frac{4-n}{2(n+2)^2},   &  \xi_p    & = \sqrt{n+2}e_p,   & \xi_{n+p}   & = -\sqrt{n+2}e_p \label{5dckf:wp}
\end{align}
where $j=1, \ldots, l$ with $l=n(n-1)/2$ and $p=1, \ldots, n$. The vector ${\mathbf 0}_{n}$ is the zero column of size $n$. The total number of the cubature points is $N=2n^2+1$. Finally, the 5D-CKF cubature vectors are generated by the following common rule
\begin{align}
 {\cal X}_{i} & = \mu + \sqrt{\Sigma}\gamma_i,  & \gamma_i & \in \left[\xi_0, \{\xi^+_j\}, \{\xi^-_j\}, \{\xi_p\}\right],  i = 0, \ldots N-1, \nonumber \\
  & \mbox{with} & w_i & \in \left[w_0, \{w^+_j\}, \{w^-_j\}, \{w_p\}\right].  \label{eq:ckf:nodes}
\end{align}

The measurement update step of the 5D-CKF estimator obeys the following equation~\cite{2013:Automatica:Jia,2018:Haykin}:
\begin{align}
\hat x_{k|k} & =\hat x_{k|k-1}+{K}_k(z_k-\hat z_{k|k-1}), \label{ckf:state}\\
{K}_{k} & =P_{xz,k}R_{e,k}^{-1}, \quad  P_{k|k} = P_{k|k-1} - {K}_k R_{e,k} {K}_k^{\top} \label{ckf:gain}
\end{align}
where the covariance matrices are calculated by
\begin{align}
\hat z_{k|k-1}&=\sum_{i=0}^{N-1} w_i  h(k,{\cal X}_{i,k|k-1}), \label{eq2.11}\\
R_{e,k}&=\sum_{i=0}^{N-1} w_i \bigl[h(k,{\cal X}_{i,k|k-1}\bigr)-\hat z_{k|k-1}\bigr]\nonumber\\
&\times \bigl[h\bigl(k,{\cal X}_{i,k|k-1}\bigr)-\hat z_{k|k-1}\bigr]^\top+R_k, \label{eq2.12}\\
P_{xz,k}&=\sum_{i=0}^{N-1} w_i \bigl[{\cal X}_{i,k|k-1}-\hat x_{k|k-1}\bigr]\bigl[h\bigl(k,{\cal X}_{i,k|k-1}\bigr)-\hat z_{k|k-1}\bigr]^\top \label{eq2.13}
\end{align}
 and the 5D-CKF cubature nodes ${\cal X}_{i,k|k-1}$, $i=0, \ldots N-1$, are defined by~\eqref{eq:ckf:nodes} with $\mu:= \hat x_{k|k-1}$ and $\sqrt{\Sigma}:=P_{k|k-1}^{1/2}$ at $t_k$, where the filter covariance matrix $P_{k|k-1}$ is factorized $P_{k|k-1}=P_{k|k-1}^{1/2}P_{k|k-1}^{{\top}/2}$.

To implement the discussed mixed-type filter effectively by using the MATLAB language, one needs to express the required formulas in a matrix-vector manner for vectorizing the calculations. For that, we define the weight column vector $w$ from the coefficients in~\eqref{eq:ckf:nodes} as well as the matrix ${\mathbb W}$ as follows:
\begin{align}
w& =\bigl[w_0,\ldots,w_{N-1}\bigr]^\top, \quad N = 2n^2+1, \label{5DCKF:w} \\
{\mathbb W} & = \bigl[I_{N}-\mathbf{1}_{N}^\top\otimes w\bigr]\mbox{\rm diag}\bigl\{w_0\ldots w_{N-1}\bigr\} \bigl[I_{N}-\mathbf{1}_{N}^\top\otimes w\bigr]^\top
\label{5DCKF:W}.
\end{align}
Additionally, we set a matrix shaped from the 5D-CKF vectors, which are located by columns ${\mathbb  X}_{k|k-1}:=\bigl[{\cal X}_{0,k|k-1}|\ldots|{\cal X}_{N-1,k|k-1} \bigr]$, i.e. it is of size $n\times N$. The symbol $\otimes$ is the Kronecker tensor product, ${\mathbf 1}_{N}$ is the unitary column of size $N$ and $I_{N}$ is the $N\times N$ identity matrix.

Following~\cite[p.~1638]{2007:Sarkka}, formulas~\eqref{eq2.11}~-- \eqref{eq2.13} can be simply written in the following matrix-vector manner
\begin{align}
{\mathbb  Z}_{k|k-1} & =h\bigl(k,{\mathbb  X}_{k|k-1}\bigr), & \hat z_{k|k-1} & = {\mathbb  Z}_{k|k-1}w, \label{ckf:zpred} \\
R_{e,k} & ={\mathbb Z}_{k|k-1}{\mathbb W}{\mathbb Z}_{k|k-1}^{\top}+R_k, &  P_{xz,k} & ={\mathbb X}_{k|k-1}{\mathbb W}{\mathbb Z}_{k|k-1}^{\top}.  \label{ckf:rek} \end{align}

In summary, to obtain the  mixed-type EKF-5DCKF estimator, one combines the continuous-discrete MATLAB-based EKF time update step with the 5D-CKF cubature rule at the measurement update. The obtained estimation method is summarized in the form of pseudo-code in Algorithm~1. Recall, it provides accurate calculations at the filters' propagation step due to the discretization error control involved in any MATLAB ODEs solver utilized, meanwhile the accurate measurement update step is due to the advanced 5D-CKF rule implemented for computing the filtered estimates.

\begin{codebox}
\Procname{{\bf Algorithm~1}. $\proc{MATLAB-based EKF-5DCKF/EKF-UKF}$}
\zi \textsc{Initialization\footnote{Without loss of generality, the estimation implementation method can be started from the filter's initial values instead of a sample draw from ${\cal N}(\bar x_0,\Pi_0)$; see the textbooks~\cite[Appendix~II.D]{1977:Bierman:book}, \cite[Theorem 9.2.1]{2000:Kailath:book}, \cite[Section~5.5]{2006:Simon:book}, \cite[Table~5.3]{2015:Grewal:book}.}:} Set $\hat x_{0|0} = \bar x_0$ and $P_{0|0} = \Pi_0$.
\zi \>Define $\gamma_i$ ($i=\overline{0,N-1}$) and $w \in {\mathbb R}^N$, ${\mathbb W} \in {\mathbb R}^{N\times N}$;
\zi \textsc{Time Update (TU)}: \Comment{\small\textsc{integrate on $[t_{k-1},t_{k}]$ with given tol.}}
\zi \>$[\hat x_{k|k-1},P_{k|k-1}]\leftarrow${\tt TU-EKF}\footnote{The pseudo-code is given in Appendix A.}$([\hat x_{k-1|k-1},P_{k-1|k-1}])$;
\zi \textsc{Measurement Update (MU)}:
\zi \>Factorize $P_{k|k-1}=P_{k|k-1}^{1/2}P_{k|k-1}^{\top/2}$;
\zi \>Generate ${\cal X}_{i,k|k-1}=\hat x_{k|k-1}+P_{k|k-1}^{1/2}\gamma_i$;
\zi \>Collect ${\mathbb  X}_{k|k-1}=\bigl[{\cal X}_{0,k|k-1},\ldots,{\cal X}_{N-1,k|k-1}\bigr]$;
\zi \>Predict ${\mathbb  Z}_{k|k-1}=h\bigl(k,{\mathbb  X}_{k|k-1}\bigr)$ and find $\hat z_{k|k-1}={\mathbb  Z}_{k|k-1}w$;
\zi \>Compute $R_{e,k}={\mathbb Z}_{k|k-1}{\mathbb W}{\mathbb Z}_{k|k-1}^{\top}+R_k$;
\zi \>Find $P_{xz,k}={\mathbb X}_{k|k-1}{\mathbb W}{\mathbb Z}_{k|k-1}^{\top}$ and ${K}_{k}=P_{xz,k}R_{e,k}^{-1}$;
\zi \>Update state $\hat x_{k|k}=\hat x_{k|k-1}+{K}_k(z_k-\hat z_{k|k-1})$;
\zi \>Update covariance $P_{k|k}=P_{k|k-1} - {K}_k R_{e,k} {K}_k^{\top}$.
\end{codebox}
%\setlinenumberplus{formuls:start}{1}

It is worth noting here that a mixed-type EKF-UKF estimator can be implemented in the same way as the above MATLAB-based EKF-5DCKF filter by using the general scheme presented in Algorithm~1 above.  Indeed, the UKF measurement update step has a similar form of equations~\eqref{ckf:state}~-- \eqref{eq2.13}, but instead of the 5D-CKF approximation formulas~\eqref{eq2.12}, \eqref{eq2.13}, one has to apply the UKF strategy as follows~\cite{2001:Wan:book,2000:Julier,2004:Julier}:
\begin{align}
R_{e,k}&=\sum_{i=0}^{2n} w_i^{(c)} \bigl[h(k,{\cal X}_{i,k|k-1}\bigr)-\hat z_{k|k-1}\bigr]\nonumber\\
&\times \bigl[h\bigl(k,{\cal X}_{i,k|k-1}\bigr)-\hat z_{k|k-1}\bigr]^\top+R_k, \label{eq2.12:UKF}\\
P_{xz,k}&=\sum_{i=0}^{2n} w_i^{(c)} \bigl[{\cal X}_{i,k|k-1}-\hat x_{k|k-1}\bigr]\bigl[h\bigl(k,{\cal X}_{i,k|k-1}\bigr)-\hat z_{k|k-1}\bigr]^\top \label{eq2.13:UKF}
\end{align}
where following the sampled-data UKF methodology, three pre-defined scalars $\alpha$, $\beta$ and $\kappa$ define $N=2n+1$ sigma points and their weights; e.g., see more details in~\cite{2001:Merwe}:
\begin{align}
w_0 & =\lambda/(n+\lambda), & w_i & = 1/(2n+2\lambda), \label{UKF_w}\\
w^{(c)}_0 & =\lambda/(n+\lambda)+1-\alpha^2+\beta, & w^{(c)}_i & = 1/(2n+2\lambda) \label{UKF_wc}
\end{align}
 with $i = \overline{1,2n}$ and $\lambda=\alpha^2(\kappa+n)-n$. The most comprehensive survey of alternative UKF parametrization variants existing in engineering literature can be found in~\cite{2015:Menegaz}. Again, the sigma vectors ${\cal X}_{i,k|k-1}$ are defined by~\eqref{eq:ckf:nodes}, i.e. with $\mu:= \hat x_{k|k-1}$ and $\sqrt{\Sigma}:=P_{k|k-1}^{1/2}$ at time instance $t_k$ where
\begin{align*}
\gamma_0 & = \mathbf{0}_n, & \gamma_i  & =  \sqrt{n+\lambda} e_i,   & \gamma_{n+i} & = -\sqrt{n+\lambda} e_i, \quad i=\overline{1,n}.
\end{align*}

Hence, one may construct a mixed-type MATLAB-based EKF-UKF estimator and implement it by using Algorithm~1 as well, but the coefficients vectors $w$ and $w^{(c)}$, the weight matrix ${\mathbb W}$ and the sigma vectors $\gamma_i$ ($\overline{0,2n}$) should be defined according to the UKF strategy. More precisely, the weight columns and matrix within the UKF estimation framework are defined as follows:
\begin{align}
w& =\bigl[w_0,\ldots,w_{N-1}\bigr]^\top \mbox{ and } w^{(c)} =\bigl[w^{(c)}_0,\ldots,w^{(c)}_{N-1}\bigr]^\top, \label{UKF:w}  \\
{\mathbb W} & = \bigl[I_{N}-\mathbf{1}_{N}^\top\otimes w\bigr]\mbox{\rm diag}\bigl\{w^{(c)}_0,\ldots, w^{(c)}_{N-1}\bigr\} \nonumber \\
& \times \bigl[I_{N}-\mathbf{1}_{N}^\top\otimes w\bigr]^\top, \quad N=2n+1. \label{UKF:W}
\end{align}
%where $N=2n+1$ is a number of sigma-points and $n$ is a dimension of the state vector to be estimated.

In summary, Algorithm~1 exposes a general continuous-discrete implementation scheme with implied MATLAB ODEs solvers for both the EKF-5CKF and EKF-UKF estimators. We stress that this common implementation methodology is of {\it covariance}-type because the filters' error covariance matrix is propagated and updated while estimation. A dual class of implementation methods exists in the engineering literature and assumes a processing of information matrix, say $\Lambda_{k|k}=P_{k|k}^{-1}$, instead of the related covariance $P_{k|k}$. The {\it information}-type filtering solves the state estimation problem without {\it a priori} information~\cite[p.~356-357]{2015:Grewal:book}. In this case the initial error covariance matrix $\Pi_0$ is too ``large'' and, hence, the initialization step $\Pi_0 := \infty$ yields various complications for covariance filtering, while the information algorithms simply imply $\Lambda_0 := 0$.
It is worth noting here that other approaches for solving the filter initialization problem include various techniques intended for extracting a prior information from available measurement data; e.g., see the single-point
or two-point filter initialization procedure from the first measurement or first two measurements discussed in~\cite[pp.~246-248]{1993:Bar-Shalom}, \cite[pp.~115-118]{2004:Ristic}, \cite[pp.35-39]{2012:Mallick:book}.

Finally, the information-form filtering is of special interest for practical applications where the number of states to be estimated is much less than the dimension of measurement vector $m>>n$. In this case the measurement update step in its traditional {\it covariance}-form formulation is extremely expensive because of the required $m$-by-$m$ matrix inversion, but the related {\it information}-form is simple~\cite[Section~7.7.3]{2015:Grewal:book}. A clear example of applications where information-form algorithms greatly benefit is the distributed multi-sensor architecture applications with
a centralized KF used for processing the information coming from a large number of sensors~\cite{2010:Sayed}. The information-type methods and their stable square-root implementations are derived for the classical KF in~\cite{1995:Park,1969:Dyer,1986:Oshman,1988:Oshman,2019:ECC:TsyganovaKulikova,2020:EJC:Tsyganova} {\it etc.}, for the UKF estimation approach in~\cite{2008:Lee:infUKF,2012:Liu:infUKF} {\it etc.}, for the 3D-CKF strategy in~\cite{2017:Arasaratnam,2013:Chandra} with a comprehensive overview given in~\cite[Chapter~6]{2019:Chandra:book}. To the best of the authors' knowledge, the information-type 5D-CKF methods still do not exist in engineering literature. Their derivation could be an area for a future research.

\section{General MATLAB-based square-rooting approach} \label{main:result}

To ensure a symmetric form and positive definiteness of the filters' error covariance in a finite precision arithmetics, the {\it square-root} methods are traditionally derived within the Cholesky decomposition; see~\cite{2000:Kailath:book,2006:Simon:book,2015:Grewal:book} and many other studies. Assume that $P_{k|k}=P_{k|k}^{1/2}P_{k|k}^{\top/2}$ where $P_{k|k}^{1/2}$ is a lower triangular matrix with positive entries on the main diagonal. The key idea of any square-root implementation is to process the square-root factors $P_{k|k}^{1/2}$ instead of the entire $P_{k|k}$. For that, the original filtering equations (see Algorithm~1) should be re-derived in terms of the Cholesky factors involved.
%This computational approach ensures a symmetric form and positive (semi-) definiteness of the filter error covariance matrix although the roundoff errors influence its square-root calculation~\cite[Chapter~7]{2015:Grewal:book}. %Taking into account the condition number of covariance matrices, Cholesky-based KF implementations are proved to maintain the computations with a double precision as shown in~\cite{1971:Kaminski}.

Unfortunately, the derivation of mathematically equivalent filtering formulas in terms of the square-root factors is not always feasible. For instance, the lack of a square-root 5D-CKF solution has been recently pointed out in~\cite{2018:Haykin}. The same problem exists for the UKF estimator. More precisely, the problem is in negative coefficients appeared in the sampled-data filtering strategies for approximating the covariance matrices. It is well known that stable orthogonal rotations are utilized for computing the Cholesky factorization of a positive definite matrix $C = A^{\top}A+B^{\top}B$, but the sampled-data estimators under examination imply equation $C = A^{\top}A \pm B^{\top}B$. Indeed, the coefficient matrix ${\mathbb W}$ contains the negative 5D-CKF weights $\{w_p\}$, $p=\overline{1,2n}$ when $n>4$ as well as the {\it classical} UKF parametrization with $\alpha=1$, $\beta=0$ and $\kappa=3-n$ yields $w^{(c)}_0<0$ when $n>3$; see~\cite{2001:Merwe}. Previously, the UKF square-rooting problem has been solved by utilizing
{\it one-rank} Cholesky update procedure wherever the negative UKF sigma-point coefficients are encountered~\cite{2001:Wan:book}. One may follow the same way with respect to the 5D-CKF estimator. However, as correctly pointed out in~\cite{2009:Haykin}, such methods are, in fact, the {\it pseudo-square-root} algorithms because they still require the Cholesky update procedure at each filtering step where the downdated matrix should be also positive definite. If this property is violated due to roundoff, then a failure of the one-rank Cholesky update yields the UKF/5D-CKF shutoff. Hence, such methods do not have the numerical robustness benefit of the {\it exact} square-root algorithms where the Cholesky decomposition is demanded for the inial covariance $\Pi_0>0$ only.
We resolve the 5D-CKF and UKF square-rooting problem by an elegant way suggested recently in~\cite{2020:SP:Kulikov,2020:IETSonar:KulikovaKulikov}.

\begin{dfn} (see~\cite{2003:Higham}) The $J$-orthogonal matrix $Q$ is defined as one, which satisfies $Q^{\top}JQ=QJQ^{\top}=J$ where $J=\mbox{\rm diag}(\pm 1)$ is a signature matrix. \end{dfn}

Clearly, if the matrix $J$ equals to identity matrix, then the $J$-orthogonal matrix $Q$ reduces to usual orthogonal one.
The $J$-orthogonal transformations are useful for finding the Cholesky factorization of a positive definite matrix $C \in {\mathbb R}^{s\times s}$ that obeys the equation
\begin{equation}
C = AA^{\top}-BB^{\top} \label{C:matrix:form}
\end{equation}
where $A \in {\mathbb R}^{s\times p}$, $(p\ge s)$ and $B \in {\mathbb R}^{s\times q}$.

If we can find a $J$-orthogonal matrix $Q$ such that
\begin{equation}
\begin{bmatrix}
A & B
\end{bmatrix} Q
=
\begin{bmatrix}
R & 0
\end{bmatrix} \label{hyperbolic_qr}
\end{equation}
with $J = \mbox{diag}\{I_p, -I_q\}$ and $R \in {\mathbb R}^{s\times s}$ is a lower triangular, then
\begin{align}
C & =
\begin{bmatrix}
A & B
\end{bmatrix}
J
\begin{bmatrix}
A & B
\end{bmatrix}^{\top}
=
\begin{bmatrix}
A & B
\end{bmatrix}
QJQ^{\top}
\begin{bmatrix}
A & B
\end{bmatrix}^{\top} \nonumber \\
&  = \Bigl[
[\underbrace{R}_s \quad \underbrace{0, \ldots, 0}_{p-s} ] \quad {\bf 0}_{s\times q}
\Bigr]
J
\Bigl[
[\underbrace{R}_s \quad \underbrace{0, \ldots, 0}_{p-s} ] \quad {\bf 0}_{s\times q}
\Bigr]^{\top} \nonumber \\
& = RI_{s}R^{\top} + {\bf 0}_{s\times (p-s)} I_{p-s} {\bf 0}_{s\times {(p-s)}}^{\top} - {\bf 0}_{s\times q} I_q {\bf 0}_{s\times q}^{\top}
= RR^{\top}, \label{proof:JQ}
\end{align}
i.e. $R$ is the lower triangular Cholesky factor, which we are looking for. It is important to satisfy $p\ge s$. The factorization in~\eqref{hyperbolic_qr} is called the {\it hyperbolic} $QR$ factorization and its effective implementation methods can be found, for instance, in~\cite{2003:Bojanczyk}.

Our first step while designing a general continuous-discrete square-root scheme, which is valid for both the UKF and/or 5D-CKF methods, is to represent the equations for computing the filters' covariance matrices (e.g., see equation~\eqref{ckf:rek}) in the required form~\eqref{C:matrix:form} and, next, to use the {\it hyperbolic} $QR$ transformation with the appropriate $J$-orthogonal matrix to get the square-root factor $R$ according to~\eqref{hyperbolic_qr}. For that, we need to factorize the coefficient matrix ${\mathbb W}$ because it is involved in the sampled-data estimators' formulas; e.g., equation~\eqref{ckf:rek} is $R_{e,k} ={\mathbb Z}_{k|k-1}{\mathbb W}{\mathbb Z}_{k|k-1}^{\top}+R_k$. The key problem at this step is the possibly negative weights and infeasible square root operation. Instead, we define the matrix $|{\mathbb W}|^{1/2}$ and the related signature matrix\footnote{For zero entries $w_{i}=0$ we set $\mbox{\rm sgn}(w_{i})=1$.} as follows:
\begin{align}
|{\mathbb W}|^{1/2}& = \bigl[I_{N}-\mathbf{1}_{N}^\top\otimes w\bigr] \mbox{\rm diag}\left\{\sqrt{|w_0|},\ldots, \sqrt{|w_{N-1}|}\right\}, \label{class_W}\\
S & = \mbox{\rm diag}\left\{\mbox{\rm sgn}(w_{0}),\ldots,\mbox{\rm sgn}(w_{N-1})\right\} \label{class_S}
\end{align}
where the same approach is utilized for the UKF estimation strategy, except that $|w^{(c)}|$ replaces the vector $|w|$ in the matrix $|{\mathbb W}|^{1/2}$ calculation; see formula~\eqref{UKF:W}.

{\bf Case 1.} According to the 5D-CKF rule, we have $2n$ negative weights for all problems of size $n>4$; see formula~\eqref{5dckf:wp}. Thus, we get $q:=2n$ in~\eqref{C:matrix:form} and transformation~\eqref{hyperbolic_qr}. The related 5D-CKF signature matrix has a form of $S = \mbox{diag}\{I_{N-2n}, -I_{2n}\}$ where $N=2n^2+1$ is the total number of the 5D-CKF points. Following equation~\eqref{eq:ckf:nodes}, the negative coefficients are all located at the end of the coefficient vector $w$. Thus, the equation for computing the covariance matrix in~\eqref{eq2.12} has the required representation~\eqref{C:matrix:form} that is suitable for computing the Cholesky factor of the matrix $R:=R_{e,k} \in {\mathbb R}^{m\times m}$. Hence, we conclude that $s:=m$ in~\eqref{C:matrix:form}, \eqref{hyperbolic_qr} and recall $q:=2n$. We can apply the hyperbolic QR as explained in~\eqref{proof:JQ} to the pre-array $A \in {\mathbb R}^{s\times p}$ with $s:=m$, $p:=m+N$ (satisfy $p\ge s$) to find the post-array $R \in {\mathbb R}^{s\times s}$, which is the required Cholesky factor $R_{e,k}^{1/2} \in {\mathbb R}^{m\times m}$, as follows:
\begin{align}
 R_{e,k} & = R_{e,k}^{1/2}R_{e,k}^{\top/2} = {\mathbb Z}_{k|k-1}{\mathbb W}{\mathbb Z}_{k|k-1}^{\top}+R_k \nonumber \\
 & = R_k^{1/2}R_k^{\top/2}+ {\mathbb Z}_{k|k-1}|{\mathbb W}|^{1/2}S|{\mathbb W}|^{\top/2}{\mathbb Z}_{k|k-1}^{\top}\nonumber  \\
 & = \left[R_k^{1/2}, {\mathbb Z}_{k|k-1}|{\mathbb W}|^{1/2}\right]\mbox{diag}\{I_m,S\} \left[R_k^{1/2}, {\mathbb Z}_{k|k-1}|{\mathbb W}|^{1/2}\right]^{\top}\nonumber \\
 & = \left[R_k^{1/2}, {\mathbb Z}_{k|k-1}|{\mathbb W}|^{1/2}\right]\mbox{diag}\{I_m,I_{N-2n},-I_{2n}\} \nonumber \\
 & \times \left[R_k^{1/2}, {\mathbb Z}_{k|k-1}|{\mathbb W}|^{1/2}\right]^{\top} = \left[R_k^{1/2}, \left[{\mathbb Z}_{k|k-1}|{\mathbb W}|^{1/2}\right]_{1:N-2n}\right]\nonumber  \\
 & \times \mbox{diag}\{I_m,I_{N-2n}\}\left[R_k^{1/2}, \left[{\mathbb Z}_{k|k-1}|{\mathbb W}|^{1/2}\right]_{1:N-2n}\right]^{\top} \nonumber \\
 &  - \left[{\mathbb Z}_{k|k-1}|{\mathbb W}|^{1/2}\right]_{N-2n+1:N}I_{2n}\left[{\mathbb Z}_{k|k-1}|{\mathbb W}|^{1/2}\right]_{N-2n+1:N}^{\top}  \label{Rek_trans}
\end{align}
where the $J$-orthogonal matrix has a form of $J:=\mbox{diag}\{I_m,S\}=\mbox{diag}\{I_m,I_{N-2n}, -I_{2n}\}$ and $\left[{\mathbb Z}_{k|k-1}|{\mathbb W}|^{1/2}\right]_{1:N-2n}$ denotes the block collected from the first $N-2n$ columns of the matrix product ${\mathbb Z}_{k|k-1}|{\mathbb W}|^{1/2}$ and $\left[{\mathbb Z}_{k|k-1}|{\mathbb W}|^{1/2}\right]_{N-2n+1:N}$ is the submatrix consisting of the last $2n$ columns of the mentioned matrix product.

{\bf Case 2.} For the UKF filters, one should pay an attention to the chosen UKF parametrization~\eqref{UKF_w}, \eqref{UKF_wc} and, in particular, to the location of negative values in $w^{(c)}$ utilized for approximating the covariance matrices under the UKF strategy; e.g., see equation~\eqref{eq2.12:UKF}. To be able to apply {\it hyperbolic} $QR$ transformations for computing the required Cholesky square-root factors, one needs to satisfy formula~\eqref{C:matrix:form}. In other words, all  negative weights in vector $w^{(c)}$ should be located at the end of $w^{(c)}$ that corresponds to the term $-BB^{\top}$ in representation~\eqref{C:matrix:form}. Thus, one needs to sort both the vector $w^{(c)}$ and the vector $w$ at the same way prior to filtering. Alternatively, the utilized hyperbolic QR method might be automatic in the sense that it should involve permutations to ensure that the post-array $R$ is indeed the required Cholesky factor, i.e. it is a lower triangular matrix with the positive entries on the diagonal; e.g., see methods in~\cite{2020:EJC:Kulikov}.
In this work, we apply the hyperbolic QR transformation from~\cite{2003:Bojanczyk} that is valid for equation~\eqref{C:matrix:form}. Thus, we do permutations prior to filtering and, hence, the coefficient matrix ${\mathbb W}$ and its $|{\mathbb W}|$ are defined as usual where the corresponding signature matrix now contains the negative entries located at the end of the main diagonal. The permutations are done at the initial filtering step and this operation does not change the summation result in the UKF formulas~\eqref{eq2.11}, \eqref{eq2.12:UKF} and \eqref{eq2.13:UKF}. For instance, the {\it classical} UKF parametrization suggested in~\cite{2001:Merwe} implies $\alpha=1$, $\beta=0$ and $\kappa=3-n$ and, hence, $w_0^{(c)}<0$ for all estimation problems of size $n>3$. Following our square-rooting solution, we locate this negative value at the end of the vectors $w^{(c)} = [w^{(c)}_{1},\ldots,w^{(c)}_{2n},w^{(c)}_{0}]$ and $w = [w_{1},\ldots,w_{2n},w_{0}]$. Thus, we get $q:=1$ in transformation~\eqref{hyperbolic_qr} and the UKF signature matrix under the classical parametrization has a form $S = \mbox{diag}\{I_{2n}, -1\}$. Next, the same derivation way~\eqref{Rek_trans} holds for computing the Cholesky factor $R_{e,k}$ by factorizing the UKF equation~\eqref{eq2.12:UKF} and utilizing the hyperbolic QR factorization with the appropriate $J$-orthogonal matrix $J:=\mbox{diag}\{I_m,S\}=\mbox{diag}\{I_m, I_{2n}, -1\}$.

Next, we design  elegant and simple {\it array} square-root formulas for updating the filters' Cholesky matrix factors involved.  Let us consider the pre-array ${\mathbb A}_k$ and block lower triangularize it as follows:
\begin{equation*}
\underbrace{
\begin{bmatrix}
R_k^{1/2}  & {\mathbb  Z}_{k|k-1}|{\mathbb W}|^{1/2}\\
\mathbf{0} & {\mathbb  X}_{k|k-1}|{\mathbb W}|^{1/2}
\end{bmatrix}
}_{\rm Pre-array  \: {\mathbb A}_k} \!\!{\mathbb Q}
\! \! =\!
\underbrace{
\begin{bmatrix}
X_1 & {\bf 0}_{m\times n} & {\bf 0}_{m\times (N-n)}\\
X_2 & X_3 & {\bf 0}_{n\times (N-n)}
\end{bmatrix}
}_{\rm Post-array  \: {\mathbb R}_k} %\label{transformation}
\end{equation*}
where ${\mathbb Q}$ is any $J$-orthogonal matrix with $J=\mbox{diag}\{I_{m}, S\}$ where $|{\mathbb W}|^{1/2}$ and $S$ are defined in~\eqref{class_W} and~\eqref{class_S}, respectively.

Following the definition of the hyperbolic $QR$ transformation and the equality ${\mathbb A}_k J {\mathbb A}_k^{\top} = {\mathbb A}_k {\mathbb Q}_k J {\mathbb Q}_k^{\top} {\mathbb A}_k^{\top}  =  {\mathbb R}_k J {\mathbb R}_k^{\top}$, we may determine the blocks $X_1$, $X_2$ and $X_3$ by factorizing the original filtering formulas. The detailed proof is available for the UKF filtering strategy in~\cite{2020:SP:Kulikov} and for the 5D-CKF rule in~\cite{2020:IETSonar:KulikovaKulikov}, respectively. Here, we briefly mention that we get $X_1:=R_{e,k}^{1/2}$, $X_3 := P^{1/2}_{k|k}$ and $X_2 := \bar P_{xz,k} = P_{xz,k}R_{e,k}^{-{\top}/2}$. Thus, the gain matrix in~\eqref{ckf:gain} is calculated by using the available post-array blocks  $\bar P_{xz,k}$ and $R_{e,k}^{1/2}$ as follows: $K_{k}=P_{xz,k}R_{e,k}^{-1} = \bar P_{xz,k}R_{e,k}^{-1/2}$.

The final step in our derivation is the MATLAB-based square-rooting solution for the time update step of the continuous-discrete EKF  approach utilized in the novel mixed-type estimators. More precisely, the EKF MDEs in~\eqref{eq2.2} should be expressed in terms of the Cholesky factor of $P(t)$, only. The following formula is proved for the lower triangular Cholesky factor $S(t)$ in~\cite{2019:Kulikov:IJRNC}:
\begin{equation}
    {dS(t)}/{dt} = S(t)\,\Phi\bigl[A(t)+A^{\top}(t)+B(t)\bigr]  \label{eq3.13}
\end{equation}
where $P(t)=S(t)S^\top(t)$ and $A(t)=S^{-1}(t)F\bigl(t,x(t)\bigr)S(t)$, $B(t)=S^{-1}(t) G(t) Q(t) G^\top(t) S^{-\top}(t)$. The mapping $\Phi(M)$ returns a lower triangular matrix defined as follows: (1) split any matrix $M$ as  $M=\bar L + D + \bar U$ where $\bar L$ and $\bar U$ are, respectively, a strictly lower and upper triangular parts of $M$, and $D$ is its main diagonal; (2) compute $\Phi(M) = \bar L + 0.5 D$.

Thus, we can formulate the first general {\it continuous-discrete} square-root filtering scheme for both the MATLAB-based EKF-5DCKF and EKF-UKF methods.

\begin{codebox}
\Procname{{\bf Algorithm~1a}. $\proc{Cholesky-based EKF-5DCKF/EKF-UKF}$}
\zi \textsc{Initialization:} Repeat the initial step from Algorithm~1.
\zi \>Cholesky dec.: $P_{0|0} = P_{0|0}^{1/2}P_{0|0}^{\top/2}$, $Q = Q^{1/2}Q^{\top/2}$;
\zi \>Define $\gamma_i$ ($i=\overline{0,N-1}$), $w$, $|{\mathbb W}|^{1/2}$ and signature $S$;
\zi \textsc{Time Update (TU)}: \Comment{\small\textsc{integrate on $[t_{k-1},t_{k}]$ with given tol.}}
\zi \>$[\hat x_{k|k-1},P^{1/2}_{k|k-1}]\leftarrow${\tt SR TU-EKF}\footnote{The pseudo-code is given in Appendix B.}$([\hat x_{k-1|k-1},P^{1/2}_{k-1|k-1}])$;
\end{codebox}
\begin{codebox}
\zi \textsc{Measurement Update (MU)}:
\zi \>Generate ${\cal X}_{i,k|k-1}=\hat x_{k|k-1}+P_{k|k-1}^{1/2}\gamma_i$;
\zi \>Collect ${\mathbb  X}_{k|k-1}=\bigl[{\cal X}_{0,k|k-1},\ldots,{\cal X}_{N-1,k|k-1}\bigr]$;
\zi \>Predict ${\mathbb  Z}_{k|k-1}=h\bigl(k,{\mathbb  X}_{k|k-1}\bigr)$ and $\hat z_{k|k-1}={\mathbb  Z}_{k|k-1}w$;
\zi \>Build pre-array ${\mathbb A}_k$ and block lower triangulates it \label{SR:Rek}
\zi \>$ \underbrace{
\begin{bmatrix}
R_k^{1/2}  & {\mathbb  Z}_{k|k-1}|{\mathbb W}|^{1/2}\\
\mathbf{0} & {\mathbb  X}_{k|k-1}|{\mathbb W}|^{1/2}
\end{bmatrix}
}_{\rm Pre-array  \: {\mathbb A}_k} {\mathbb Q}
=
\underbrace{
\begin{bmatrix}
R_{e,k}^{1/2} & {\bf 0}_{m\times n} & {\bf 0}_{m\times (N-n)}\\
\bar P_{xz,k} & P^{1/2}_{k|k} & {\bf 0}_{n\times (N-n)}
\end{bmatrix}
}_{\rm Post-array  \: {\mathbb R}_k}$
\zi \>with ${\mathbb Q}$ is any $J=\mbox{diag}\{I_m,S\}$-orthogonal rotation;
\zi \>Extract $R_{e,k}^{1/2}$, $\bar P_{xz,k}$, $P^{1/2}_{k|k}$ and find ${K}_{k}=\bar P_{xz,k}R_{e,k}^{-1/2}$;
\zi \>Update state estimate $\hat x_{k|k}=\hat x_{k|k-1}+{K}_k(z_k-\hat z_{k|k-1})$.
\end{codebox}
%\setlinenumberplus{formuls:start}{1}

Finally, one more Cholesky-based implementation method can be designed by factorizing the symmetric Joseph-type stabilized equation derived for the sampled-data UKF rule in~\cite{2020:SP:Kulikov}:
\[
 P_{k|k}  = \left[{\mathbb X}_{k|k-1}-{K}_{k}{\mathbb Z}_{k|k-1}\right]{\mathbb W}\left[{\mathbb X}_{k|k-1}-{K}_{k}{\mathbb Z}_{k|k-1}\right]^{\top} + {K}_{k}R_k {K}_{k}^{\top}.
\]

The symmetric form of the formula above permits the following factorization:
\[
\underbrace{
\begin{bmatrix}
{K}_{k} R_k^{1/2} & \left[{\mathbb X}_{k|k-1}-{K}_{k}{\mathbb Z}_{k|k-1}\right]|{\mathbb W}|^{1/2}
\end{bmatrix}
}_{\rm Pre-array  \: {\mathbb A}_k} {\mathbb Q}=
\underbrace{
\begin{bmatrix}
P_{k|k}^{1/2} &  0
\end{bmatrix}
}_{\rm Post-array  \: {\mathbb R}_k}\]
where ${\mathbb Q}$ is any $J$-orthogonal transformation with $J=\mbox{diag}\{I_{m}, S\}$ that block triangularizes the pre-array.

Thus, an alternative MATLAB-based square-root implementation for the examined mixed-type filtering strategy can be summarized as follows.

\begin{codebox}
\Procname{{\bf Algorithm~1b}. $\proc{Cholesky-based EKF-5DCKF/EKF-UKF}$}
\zi \textsc{Initialization:} Repeat from Algorithm~1a;
\zi \textsc{Time Update (TU)}: Repeat from Algorithm~1a;
\zi \textsc{Measurement Update (MU)}:
\zi \>Generate ${\cal X}_{i,k|k-1}=\hat x_{k|k-1}+P_{k|k-1}^{1/2}\gamma_i$;
\zi \>Collect ${\mathbb  X}_{k|k-1}=\bigl[{\cal X}_{0,k|k-1},\ldots,{\cal X}_{N-1,k|k-1}\bigr]$;
\zi \>Predict ${\mathbb  Z}_{k|k-1}=h\bigl(k,{\mathbb  X}_{k|k-1}\bigr)$ and $\hat z_{k|k-1}={\mathbb  Z}_{k|k-1}w$;
\zi \>Build pre-array ${\mathbb A}_k$ and block lower triangulates it
\zi \>$\underbrace{
\begin{bmatrix}
R_k^{1/2}  & {\mathbb  Z}_{k|k-1}|{\mathbb W}|^{1/2}
\end{bmatrix}
}_{\rm Pre-array  \: {\mathbb A}_k} {\mathbb Q}
 =
\underbrace{
\begin{bmatrix}
R_{e,k}^{1/2} & {\bf 0}_{m\times N}
\end{bmatrix}
}_{\rm Post-array  \: {\mathbb R}_k}$
\zi \>with ${\mathbb Q}$ is any $J=\mbox{diag}\{I_m,S\}$-orthogonal rotation;
\zi \>Extract $R_{e,k}^{1/2}$ and find $\bar P_{xz,k}={\mathbb X}_{k|k-1}{\mathbb W}{\mathbb Z}_{k|k-1}^{\top}R_{e,k}^{-{\top}/2}$;
\zi \>Calculate the filter gain ${K}_{k}=\bar P_{xz,k}R_{e,k}^{-1/2}$;
\zi \>Build pre-array ${\mathbb A}_k$ and block lower triangulates it \label{SR:Rek}
\zi \>$ \underbrace{
\begin{bmatrix}
{K}_{k} R_k^{1/2} & \left[{\mathbb X}_{k|k-1}-{K}_{k}{\mathbb Z}_{k|k-1}\right]|{\mathbb W}|^{1/2}
\end{bmatrix}
}_{\rm Pre-array  \: {\mathbb A}_k} {\mathbb Q}=
\underbrace{
\begin{bmatrix}
P_{k|k}^{1/2} &  0
\end{bmatrix}
}_{\rm Post-array  \: {\mathbb R}_k}$
\zi \>with ${\mathbb Q}$ is any $J=\mbox{diag}\{I_m,S\}$-orthogonal rotation;
\zi \>Extract $P_{k|k}^{1/2}$ and find $\hat x_{k|k}=\hat x_{k|k-1}+{K}_k(z_k-\hat z_{k|k-1})$.
\end{codebox}

\section{Numerical experiments} \label{numerical:experiments}

For providing a fair comparative study of the novel MATLAB-based methods with existing nonlinear KF-like methods and the previously published NIRK-based filters in~\cite{2020:SP:Kulikov,2020:IETSonar:KulikovaKulikov}, we examine the same test problem, which was suggested for illustrating a performance of various CKF techniques in~\cite{2010:Haykin,2018:Haykin}. Additionally, we examine both the Gaussian noise scenario and non-Gaussian one presented in~\cite{KuKu18SP,KuKu20aIFAC}. In the cited papers, it is shown that the UKF and CKF filtering approaches provide a high estimation quality and, as a result, both techniques are applicable for processing the nonlinear estimation problems with non-Gaussian uncertainties as well.

\begin{exampe} \label{ex:1} When performing a coordinated turn in the horizontal plane, the aircraft's dynamics obeys equation~\eqref{eq1.1} with the following drift function and diffusion matrix
\[
f(\cdot)=\left[\dot{\epsilon}, -\omega \dot{\eta}, \dot{\eta}, \omega \dot{\epsilon}, \dot{\zeta},  0, 0\right],
G={\rm diag}\left[0,\sigma_1,0,\sigma_1,0,\sigma_1,\sigma_2\right]
\] where $\sigma_1=\sqrt{0.2}\mbox{ \rm m/s}$, $\sigma_2=0.007^\circ/\mbox{\rm s}$ and $\beta(t)$ is the {\it standard} Brownian motion, i.e. $Q=I_7$.
The state vector consists of seven entries, i.e. $x(t)= [\epsilon, \dot{\epsilon}, \eta, \dot{\eta}, \zeta, \dot{\zeta}, \omega]^{\top}$, where $\epsilon$, $\eta$, $\zeta$ and $\dot{\epsilon}$, $\dot{\eta}$, $\dot{\zeta}$ stand for positions and corresponding velocities in the Cartesian coordinates at time $t$, and $\omega(t)$ is the (nearly) constant turn rate. The initial conditions are $\bar x_0=[1000\,\mbox{\rm m}, 0\,\mbox{\rm m/s}, 2650\,\mbox{\rm m},150\,\mbox{\rm m/s}, 200\,\mbox{\rm m}, 0\,\mbox{\rm m/s},\omega^\circ/\mbox{\rm s}]^{\top}$ and $\Pi_0=\mbox{\rm diag}(0.01\,I_7)$. We fix the turn rate to $\omega=3^\circ/\mbox{\rm s}$.

{\bf Case 1: The original problem with Gaussian noise}. The measurement model is taken from~\cite{2010:Haykin}, but accommodated to MATLAB as follows:
\[
\begin{bmatrix}
r_k \\
\theta_k \\
\phi_k
  \end{bmatrix}
 =
  \begin{bmatrix}
  \sqrt{\epsilon^2_k+\eta^2_k+\zeta^2_k} \\
  {\rm atan2}\left({\eta_k},{\epsilon_k}\right) \\
  {\rm atan}\left(\frac{\zeta_k}{{\sqrt{\epsilon^2_k+\eta^2_k}}}\right)
  \end{bmatrix}
  + v_k,
  \begin{array}{l}
    v_k \sim {\cal N}(0,R); \\
    R  ={\rm diag}(\sigma_r^2,\sigma_\theta^2,\sigma_\phi^2)
\end{array}
\]
where the implementation of the MATLAB command {\tt atan2} in target tracking is explained in~\cite{2007:Brehard} in more details.
The observations $z_k = [r_k, \theta_k, \phi_k]^{\top}$, $k=1,\ldots, K$ come at some constant sampling intervals $\Delta = [t_{k-1}, t_k]$. The radar is located at the origin and is equipped to measure the range $r_k$, azimuth angle~$\theta$ and the elevation angle~$\phi$. The measurement noise covariance matrix is constant over time and $\sigma_r=50$~m, $\sigma_\theta=0.1^\circ$, $\sigma_\phi=0.1^\circ$.

{\bf Case 2: The original problem with Glint noise scenario}. Following \cite{HeMa87,Wu93,BiTa06}, the glint noise measurement uncertainty is modeled by the sum of two Gaussian random variables
\begin{equation}\label{eq5.1}
v_k=(1-p_g){\mathcal N}(0,R)+p_g{\mathcal N}(0,R_g)
\end{equation}
where the constant $p_g$ stands for the probability of the glint, the matrix $R$ refers to the measurement noise covariance in the above Gaussian noise target tracking scenario and $R_g$ denotes the covariance of the glint outliers. In line with the cited literature, we cover the scenario with the glint probability $p_g:=0.25$ and the glint noise covariance $R_g=100 R$.
\end{exampe}

\begin{table}[ht!]
\renewcommand{\arraystretch}{1.3}
\caption{The $\mbox{\rm ARMSE}_p$ (m) of various filtering methods in Example~\ref{ex:1}, Case~1.} \label{tab:acc}
{\scriptsize
\begin{tabular}{r|r|r|r|r|r}
\hline
$\Delta$(s) & \multicolumn{3}{c|}{\bf Existing methods} & \multicolumn{2}{c}{\bf New mixed-type methods}  \\
\cline{2-6}
& 5D-CKF & UKF; see  & EKF; see & EKF-UKF & EKF-5DCKF  \\
& in~\cite{2018:Haykin} & \cite[Alg.~1]{2020:EJC:Kulikov}  & \cite[Alg.~1]{KuKu18cIEEE_ICSTCC} &  &   \\
\hline
    1 &    62.75	 &    62.76 	 &    75.03 	 &    71.33 	 &    71.32 \\ 	
    2 &    83.34	 &    83.39 	 &   116.40 	 &    99.69 	 &    99.69 \\ 	
    3 &    90.91	 &    90.97 	 &   171.00 	 &   108.61 	 &   108.60 \\ 	
    4 &    95.48	 &    95.59 	 &   225.90 	 &   120.60 	 &   120.60 \\	
    5 &    96.78	 &    96.82 	 &   254.80 	 &   119.20 	 &   119.20 \\	
    6 &   110.50	 &   110.50 	 &   284.10 	 &   137.73 	 &   137.60 \\	
    7 &   102.40	 &   102.50 	 &   284.50 	 &   127.50 	 &   127.50 \\	
    8 &   122.60	 &   122.90 	 &   298.80 	 &   148.31 	 &   148.30 \\	
    9 &   136.30	 &   135.20 	 &   301.50 	 &   153.30 	 &   153.30 \\	
   \hline
   10 & --- {\bf fails}  &   154.40 	         & $>$500 {\bf fails} 	 &   154.30 	 &   154.30 \\	
   11 &  	             &   184.40 	         &  	 &   157.60 	 &   157.60\\ 	
   12 & 	             &  $>$500 {\bf fails} 	 &  	 &   170.40 	 &   170.50 \\	
   \hline
\end{tabular}
}
\end{table}

In our first set of numerical experiments, we solve the original state estimation problem from~\cite{2010:Haykin} (see Case~1 study in Example~1) on the interval $[0s, 150s]$ with various sampling periods $\Delta =1, 2, \ldots, 12(s)$. The Euler-Maruyama method with the small step size $0.0005 (s)$ is utilized for computing the exact trajectory, i.e. $x^{true}_k$, and for simulating the related ``true'' measurements at the same time instants. Next, the inverse (filtering) problem is solved by various filtering methods in order to obtain the estimates $\hat x_{k|k}$, $k=1, \ldots, K$. Additionally, we compute the accumulated root mean square error (ARMSE) by averaging over $M=100$ Monte
Carlo runs. Following~\cite{2010:Haykin}, the ARMSE in position ($\mbox{\rm ARMSE}_p$) is calculated and the filters' failure is defined when $\mbox{\rm ARMSE}_p> 500$ (m).

First of all, we examine the mixed-type EKF-UKF and EKF-5DCKF estimators against the following traditional filtering methods: (i) the EKF estimator; see the algorithm summarized in~\cite[Section~II-A, Algorithm~1]{KuKu18cIEEE_ICSTCC}, (ii) the UKF method; see the details in~\cite[Section~II-B]{KuKu16IEEE_TSP} and also its brief algorithmic representation in~\cite[Section 2, Algorithm~1]{2020:EJC:Kulikov}, and (iii)  the 5D-CKF suggested in~\cite{2018:Haykin}. These previously published {\it continuous-discrete} filters are fixed step size methods based on the It$\hat{\rm o}$-Taylor expansion of order 1.5, as suggested in~\cite{2010:Haykin}. They require a number of subdivisions for each sampling interval to be supplied by the user prior to filtering. We implement these estimators with $L=64$ subdivisions. The newly-suggested MATLAB-based EKF-UKF and EKF-5DCKF algorithms do not require the number of subdivisions but they demand the tolerance value to be given prior to filtering. We implement the novel mixed-type methods with $\epsilon_g = 10^{-4}$ and with the use of the MATLAB ODEs solver \verb"ode45". Finally, all estimators under examination are tested with the same initial conditions, with the same simulated ``true'' state trajectory and the same measurement data. The results of this set of numerical experiments are summarized in Table~\ref{tab:acc} where the estimation accuracies are collected.

Having analyzed the obtained results, we make a few conclusions. Firstly, we observe that the EKF method is the less accurate estimator among all tested methods. Besides, it fails to solve the stated estimation problem for a long sampling intervals. Indeed, the computed $\mbox{\rm ARMSE}_p$ is more than $500$ (m) for sampling interval $\Delta > 10$ (s) that means its failure as indicated in~\cite{2010:Haykin,2018:Haykin}. Secondly, the estimation quality of the 5D-CKF estimator is very similar to the UKF accuracies but the 5D-CKF strategy slightly outperforms the UKF framework. However, our numerical study illustrates that the 5D-CKF estimator is less stable than the UKF algorithm. It is clearly seen from the sign `---' appeared for $\Delta = 10$ (s). This means that the 5D-CKF method fails because of the round-off errors and due to unfeasible Cholesky decomposition required for generating the 5D-CKF cubature vectors. This result is in line with the numerical study presented in~\cite{2018:Haykin} where the instability issue of the 5D-CKF estimator is emphasized. In particular, this observation evidences the theoretical rigor of the square-rooting procedure for the 5D-CKF methods. The importance of using square-root approach while implementing the 5D-CKF is stressed in~\cite{2018:Haykin}.

\begin{table*}[ht!]
\renewcommand{\arraystretch}{1.3}
\caption{The performance of various mixed-type EKF-UKF implementation methods in Example~\ref{ex:1}.} \label{Tab:1}
\centering
{\small
\begin{tabular}{r||c|c||c|c||c|c}
\hline
&  \multicolumn{2}{c||}{Conventional EKF-UKF} &  \multicolumn{4}{c}{Square-root EKF-UKF implementations} \\
\cline{2-7}
 & NIRK-based  &  MATLAB-based & NIRK-based & MATLAB-based  & NIRK-based &  MATLAB-based\\
 & in~\cite[Alg.~1]{2020:SP:Kulikov}  & new Alg.~1 & in~\cite[Alg.~2]{2020:SP:Kulikov} & new Alg.~1a & in~\cite[Alg.~3]{2020:SP:Kulikov}   & new Alg.~1b\\
\hline
\hline
 &  \multicolumn{6}{l}{{\bf Case~1 study: Gaussian uncertainties}. The utilized MATLAB ODEs solver is \texttt{ode45}} \\
\cline{2-7}
$\mbox{\rm ARMSE}_{p} (\mbox{\rm m})$ & 7.0030e+01	 & 7.1320e+01	 & 7.0030e+01	 & 7.1360e+01	 & 7.0030e+01	 & 7.1350e+01	\\
$\mbox{\rm ARMSE}_{v} (\mbox{\rm m/s})$ & 1.4000e+02	 & 1.4650e+02	 & 1.4000e+02	 & 1.4670e+02	 & 1.4000e+02	 & 1.4670e+02 \\	
$\mbox{\rm CPU Time} (\mbox{\rm s})$ & 1.0800e+00	 & 5.5750e-01	 & 1.3550e+00	 & 1.1520e+00	 & 1.3930e+00	 & 1.1900e+00	\\
$\mbox{\rm Time Benefit} (\mbox{\rm \%})$  & --- 	 	 &    93.72\%  	  & --- 	 	 &    17.62\%  	 & --- 	 	 &    17.06\%  	\\
\hline
 &  \multicolumn{6}{l}{{\bf Case~2 study: Glint noise case}. The utilized MATLAB ODEs solver is \texttt{ode45}} \\
\cline{2-7}
$\mbox{\rm ARMSE}_{p} (\mbox{\rm m})$ & 1.2480e+02	 & 1.3000e+02	 & 1.2480e+02	 & 1.3020e+02	 & 1.2480e+02	 & 1.3010e+02 \\	
$\mbox{\rm ARMSE}_{v} (\mbox{\rm m/s})$ & 2.3370e+02	 & 2.6080e+02	 & 2.3370e+02	 & 2.6120e+02	 & 2.3370e+02	 & 2.6110e+02 \\	
$\mbox{\rm CPU Time} (\mbox{\rm s})$ & 1.0330e+00	 & 5.6750e-01	 & 1.3030e+00	 & 1.2500e+00	 & 1.3360e+00	 & 1.2690e+00	\\
$\mbox{\rm Time Benefit} (\mbox{\rm \%})$  & --- 	 	 &    82.03\%  	  & --- 	 	 &     4.24\%  	 & --- 	 	 &     5.28\%  \\
\hline
\end{tabular}
}
\end{table*}

Finally, let us discuss the performance of the mixed-type estimators in comparison with the existing filters under examination. We observe that the EKF approach is less accurate and less stable than the novel mixed-type EKF-5DCKF and EKF-UKF estimators. At the same time, the 5D-CKF and UKF estimators outperform the related mixed-type counterparts for the estimation quality for short sampling intervals. In average, the 5D-CKF and UKF methods on $\approx 2-10$\% more accurate than the mixed-type filters EKF-5DCKF and EKF-CKF, but this holds true only for $\Delta <10$ (s). As can be seen, for $\Delta =10$(s) the situation is dramatically changed and the mixed-type filters start to outperform the 5D-CKF and UKF estimators for accuracy and stability. Indeed, the 5D-CKF fails for $\Delta \ge 10$(s) and the UKF estimator fails for $\Delta \ge 12$(s), meanwhile the novel MATLAB-based mixed-type filters work accurately and easily manage the stated estimation problems providing a good estimation quality regardless the sampling interval length under examination. The source of the failure is the accumulated discretization error, which ultimately destroys the entire filtering schemes for existing UKF and 5D-CKF frameworks. In contrast, the filters based on numerical solution of the related MDEs by using MATLAB ODEs solvers with built-in controlling techniques suggest an accurate estimation way for any sampling interval length, including the case of irregular sampling periods due to missing measurements, until the built-in discretization
error control involved in the solver keeps the discretization error reasonably small. 
 The adaptive nature of the numerical schemes involved  makes them flexible and convenient for using in practice. No preliminary tuning is required by users, except the tolerance value to be given. Indeed, the built-in advanced numerical integration scheme with the adaptive error control strategy ensures (in automatic mode) that the occurred discretization error is insignificant while propagating the mean and the filters' error covariance matrix. To design such accurate and stable MATLAB-based 5D-CKF estimator, the 5D-CKF MDEs should be derived, first. To the best of authors' knowledge, this is still an open problem.

\begin{table*}[ht!]
\renewcommand{\arraystretch}{1.3}
\caption{The performance of various mixed-type EKF-5DCKF implementation methods in Example~\ref{ex:1}.} \label{Tab:2}
\centering
{\small
\begin{tabular}{r||c|c||c|c||c|c}
\hline
&  \multicolumn{2}{c||}{Conventional EKF-5DCKF} &  \multicolumn{4}{c}{Square-root EKF-5DCKF implementations} \\
\cline{2-7}
 & NIRK-based  &  MATLAB-based & NIRK-based & MATLAB-based  & NIRK-based &  MATLAB-based\\
 & in~\cite[Alg.~1]{2020:IETSonar:KulikovaKulikov}  & new Alg.~1 & in~\cite[Alg.~1a]{2020:IETSonar:KulikovaKulikov} & new Alg.~1a & in~\cite[Alg.~1b]{2020:IETSonar:KulikovaKulikov}   & new Alg.~1b\\
\hline
\hline
 &  \multicolumn{6}{l}{{\bf Case~1 study: Gaussian uncertainties}. The utilized MATLAB ODEs solver is \texttt{ode45}} \\
\cline{2-7}
$\mbox{\rm ARMSE}_{p} (\mbox{\rm m})$   & 7.0030e+01	 & 7.1320e+01	 & 7.0030e+01	 & 7.1350e+01	 & 7.0030e+01	 & 7.1350e+01 \\	
$\mbox{\rm ARMSE}_{v} (\mbox{\rm m/s})$ & 1.4000e+02	 & 1.4650e+02	 & 1.4000e+02	 & 1.4670e+02	 & 1.4000e+02	 & 1.4670e+02\\ 	
$\mbox{\rm CPU Time}  (\mbox{\rm s})$   & 1.6510e+00	 & 6.5080e-01	 & 3.7480e+00	 & 3.0020e+00	 & 3.7460e+00	 & 3.0390e+00 \\	
$\mbox{\rm Time Benefit}  (\mbox{\rm \%})$ & --- 	 	 &   153.69\% 	  & --- 	 	 &    24.85\% 	 & --- 	 	 &    23.26\% \\
\hline
 &  \multicolumn{6}{l}{{\bf Case~2 study: Glint noise case}. The utilized MATLAB ODEs solver is \texttt{ode45}} \\
\cline{2-7}
$\mbox{\rm ARMSE}_{p} (\mbox{\rm m})$   & 1.2480e+02	 & 1.3000e+02	 & 1.2480e+02	 & 1.3010e+02	 & 1.2480e+02	 & 1.3010e+02	\\
$\mbox{\rm ARMSE}_{v} (\mbox{\rm m/s})$ & 2.3370e+02	 & 2.6080e+02	 & 2.3370e+02	 & 2.6110e+02	 & 2.3370e+02	 & 2.6110e+02	\\
$\mbox{\rm CPU Time}  (\mbox{\rm s})$   & 1.6970e+00     & 6.6020e-01	 & 3.7700e+00	 & 3.1310e+00	 & 3.7570e+00	 & 3.1510e+00	\\
$\mbox{\rm Time Benefit}  (\mbox{\rm \%})$ & --- 	 	 &   157.04\% 	  & --- 	 	 &    20.41\% 	 & --- 	 	 &    19.23\% \\
\hline
\end{tabular}
}
\end{table*}

Our second set of numerical experiments is focused on examining the previously published NIRK-based mixed-type filters from~\cite{2020:SP:Kulikov,2020:IETSonar:KulikovaKulikov} against their novel MATLAB-based counterparts derived in this paper. Again, all algorithms are tested at the same initial conditions, the same measurement data and the same tolerance value $\epsilon_g = 10^{-4}$ given by user. The MATLAB ODEs solver \verb"ode45" is utilized in all our novel algorithms. We stress that any other numerical integration scheme might be easily utilized instead of the \verb"ode45" function. We fix the sampling interval to $\Delta = 1(s)$ and calculate the following ARMSE by averaging over $100$ Monte Carlo runs~\cite{2010:Haykin}: i) the ARMSE in position ($\mbox{\rm ARMSE}_p$), and ii) the ARMSE in velocity ($\mbox{\rm ARMSE}_v$). We also provide the average CPU
time (in sec.)  for each filtering scheme under examination and, next, we calculate the computational benefit (in \%) of using
the MATLAB-based filtering method instead of its NIRK-based counterpart. The results obtained for the mixed-type EKF-UKF methods are summarized in Table~\ref{Tab:1} meanwhile Table~\ref{Tab:2} contains the outcomes for the mixed-type EKF-5DCKF filters.

Having analyzed the results presented in Tables~\ref{Tab:1} and~\ref{Tab:2}, we make a few important conclusions. First, we note that the novel MATLAB-based filters provide a similar estimation accuracy as the previously published NIRK-based estimators, but the novel algorithms are easy to implement and, more importantly, they are much faster compared to their NIRK-based counterparts. Indeed, the computed accumulated errors $\mbox{\rm ARMSE}_{p}$ and $\mbox{\rm ARMSE}_{v}$ are similar under the MATLAB-based and the related NIRK-based filtering approaches. The difference observed is due to the discretization error control applied by the adaptive numerical schemes at the time update step. In the previously published NIRK-based algorithms, the combined {\it local-global} discretization error control is utilized; see~\cite{2013:Kulikov:IMA}. It yields slightly more accurate estimates than those provided by the MATLAB solvers with their built-in {\it local} error control involved. However, the price to be paid is a computationally heavier filtering procedure with the increased CPU time. In fact, from Table~\ref{Tab:2} we observe that the new conventional MATLAB-based EKF-5DCKF filter is $\approx 1.5$ times faster on average than the related NIRK-based EKF-5DCKF variant suggested recently in~\cite{2020:IETSonar:KulikovaKulikov}. This conclusion holds true regardless the Gaussian and non-Gaussian estimation case study. Meanwhile, following the results summarized in Table~\ref{Tab:1}, the new conventional MATLAB-based EKF-UKF estimator is $\approx 94\%$ and $\approx 82\%$ faster than the previously proposed NIRK-based EKF-UKF counterpart for the Gaussian case and the glint noise scenario, respectively.

In engineering literature, it is well known that the square-root algorithms are slower than the related conventional (original) implementations, because of the QR factorizations involved. Hence, the CPU time benefits (in \%) of using the MATLAB-based {\it square-root} methods instead of their NIRK-based {\it square-root} variants are less than for the conventional implementations discussed above. Following the results in Table~\ref{Tab:1},  the square-root MATLAB-based EKF-UKF estimators are about $\approx 17\%$ faster compared to the square-root NIRK-based EKF-UKF methods for the Gaussian case study. Meanwhile the square-root MATLAB-based EKF-5DCKF filters are about $\approx 25\%$ and $\approx 20\%$ faster on average compared to their square-root NIRK-based counterparts in the Gaussian case and glint noise scenario, respectively; see Table~\ref{Tab:2}.

To investigate a difference in the filters' numerical robustness with respect to roundoff errors, we follow the ill-conditioned test problem proposed in~\cite[Example~7.2]{2015:Grewal:book} and discussed at the first time in~\cite[Examples 7.1 and 7.2]{1969:Dyer}. In the cited papers, a specially designed measurement scheme is suggested to be utilized for provoking the filters' numerical instability due to roundoff and for observing their divergence when the problem ill-conditioning increases.

\begin{exampe} \label{ex:2}
The dynamic state in Example~\ref{ex:1} is observed through the following measurement scheme
\begin{align*}
z_k & =
\begin{bmatrix}
1 & 1 & 1 & 1 & 1 &  1 &  1\\
1 & 1 & 1 & 1 & 1 &  1 &  1 +\delta
\end{bmatrix}
x_k +
\begin{bmatrix}
v_k^1 \\
v_k^2
\end{bmatrix}, \; R_k=\delta^{2}I_2
\end{align*}
where parameter $\delta$ is used for simulating roundoff effect. This increasingly ill-conditioned target tracking scenario assumes that $\delta\to 0$ and, hence, the residual covariance matrix $R_{e,k}$ to be inverted becomes ill-conditioned and close to zero that is a source of the numerical instability.
\end{exampe}

\begin{figure}[t!]
\includegraphics[width=0.5\textwidth]{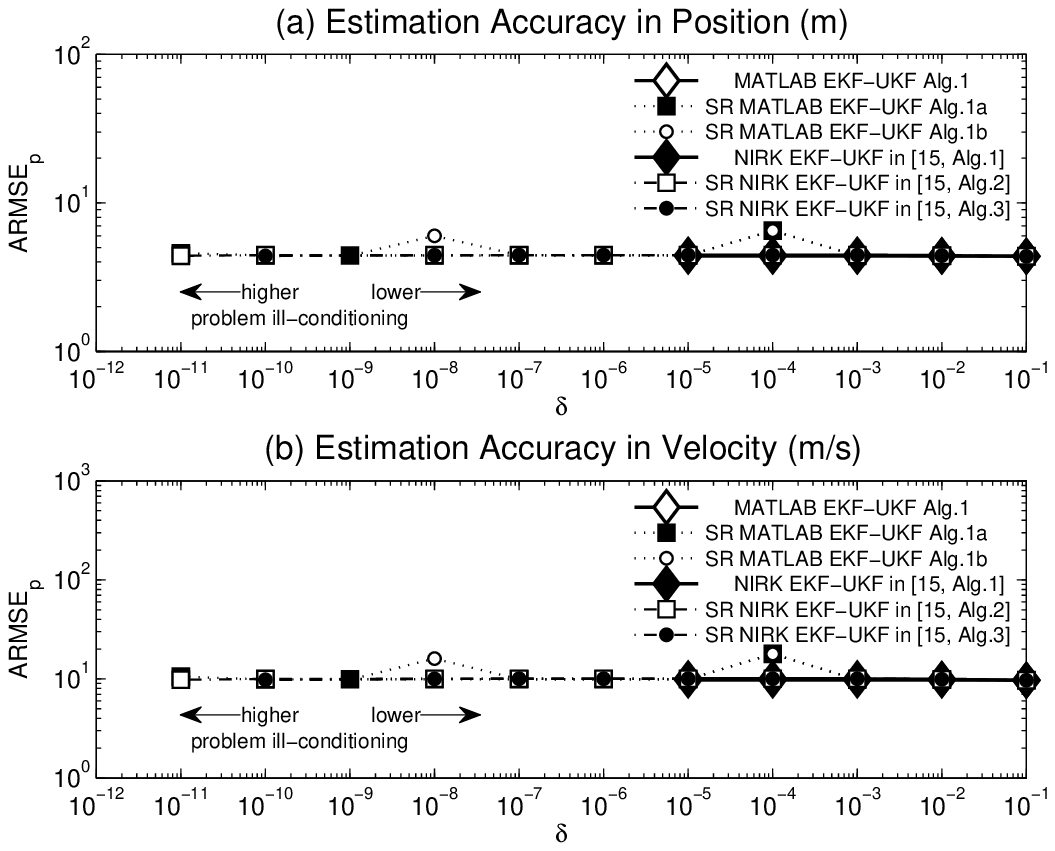}\\
\includegraphics[width=0.5\textwidth]{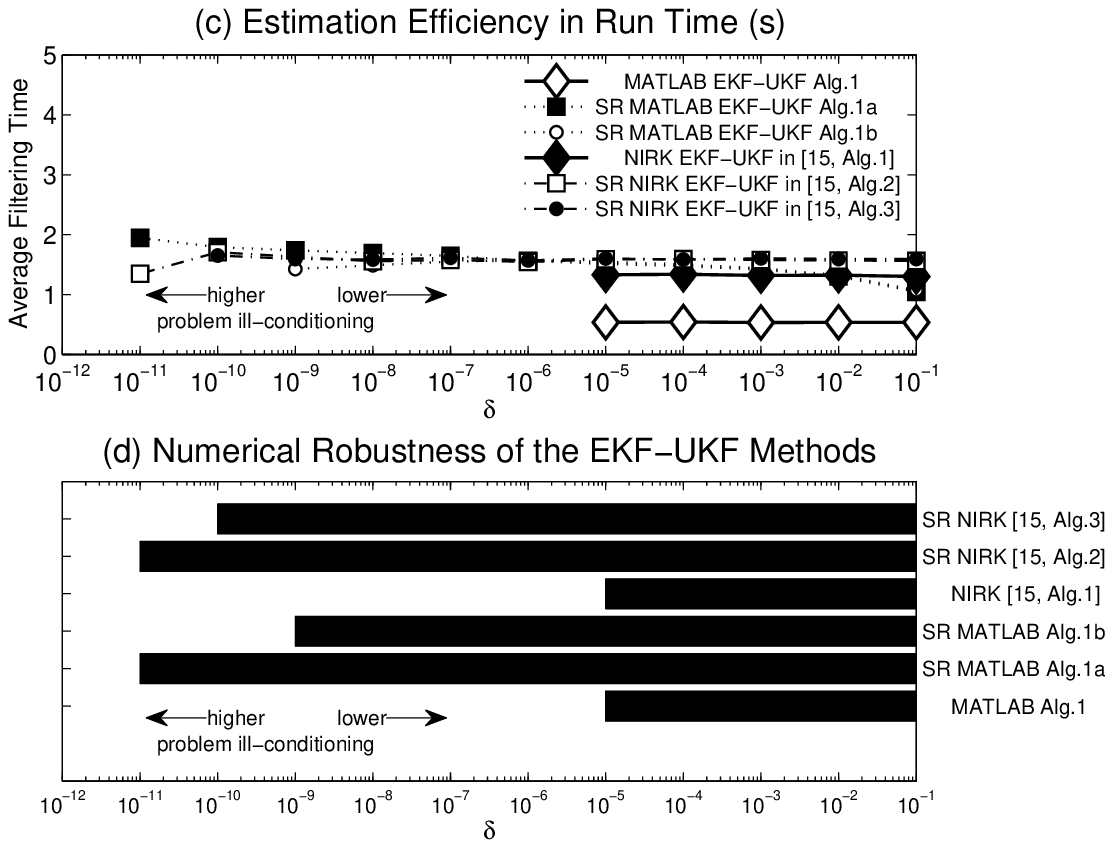}
\caption{Performance of various mixed-type continuous-discrete EKF-UKF estimators on the set of ill-conditioned test problems from Example~\ref{ex:2}.} \label{fig:1:new}
\end{figure}

Having repeated the experiment explained above for each $\delta=10^{-1},10^{-2},\ldots $, we illustrate the accumulated root mean square errors computed and CPU time (s) against the ill-conditioning parameter $\delta$ as well as the filters' breakdown value by Figs.~1 and~2. Let us consider the results obtained for the mixed-type continuous-discrete EKF-UKF filters illustrated by Fig.~\ref{fig:1:new}. We observe a similar patten in both $\mbox{\rm ARMSE}_p$ and $\mbox{\rm ARMSE}_v$ estimation accuracies degradations; see Fig.~1(a) and  Fig.~1(b). As can be seen, the conventional algorithms, i.e. both the previously published conventional NIRK-based EKF-UKF filter and the novel MATLAB-based EKF-UKF in Algorithm~1, fail at $\delta < 10^{-5}$. In other words, their numerical robustness to roundoff is the same. It is also clearly seen from Fig.~\ref{fig:1:new}(d) where the breakdown instance of each algorithm under examination is indicated. Meanwhile, from Fig.~\ref{fig:1:new}(c) it is evident that the new MATLAB-based EKF-UKF in Algorithm~1 is faster than the previously published its NIRK-based counterpart.

Next, we analyze the performance of the square-root EKF-UKF implementations. They manage to solve the ill-conditioned state estimation problems until $\delta = 10^{-11}$. At that point they still produce quite accurate estimates because the accumulated errors are small and, ultimately, they diverge at $\delta < 10^{-11}$. In summary, the new square-root MATLAB-based EKF-UKF method in Algorithm~1a and its NIRK-based variant proposed in~\cite[Alg.~2]{2020:SP:Kulikov} are equally robust with respect to roundoff errors, i.e. they diverge with the same speed. It is interesting to note that for our test problem, the MATLAB-based EKF-UKF summarized in Algorithm~1b is the fastest method to diverge among all square-root EKF-UKF implementations under examination. The source of its worse numerical robustness might be in a large number of arithmetic operation involved due to extra $QR$ factorization at the measurement update step compared to Algorithm~1a.

\begin{figure}[t!]
\includegraphics[width=0.5\textwidth]{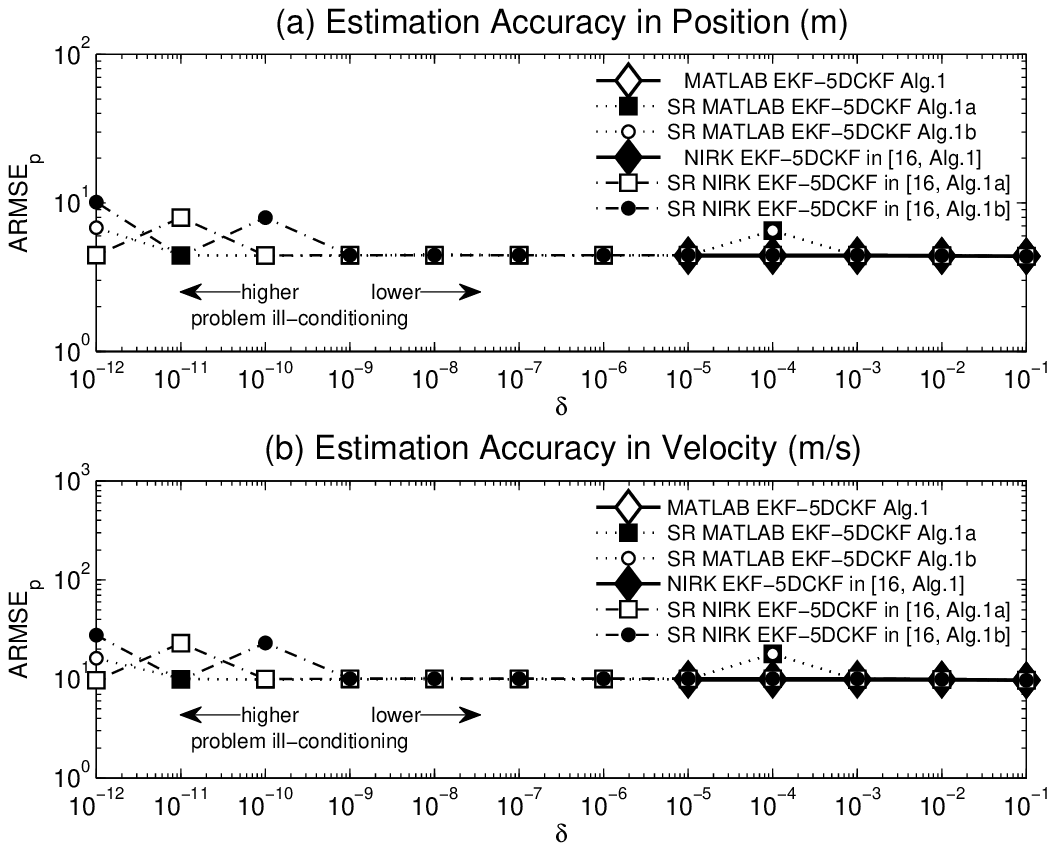}\\
\includegraphics[width=0.5\textwidth]{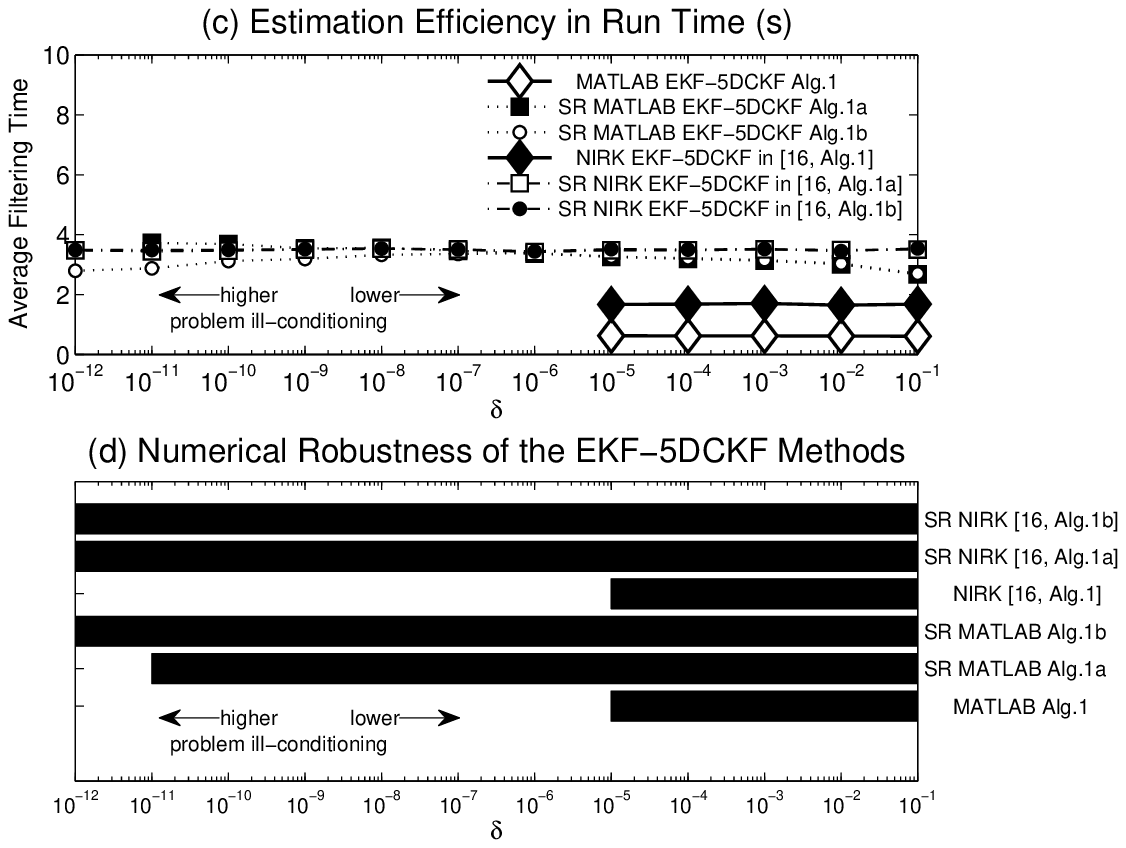}
\caption{Performance of various mixed-type continuous-discrete EKF-5DCKF estimators on the set of ill-conditioned test problems from Example~\ref{ex:2}.} \label{fig:2:new}
\end{figure}

Finally, let us consider the results obtained for the mixed-type EKF-5DCKF estimators illustrated by Fig.~2. We make the following conclusions. First, the conventional EKF-5DCKF implementation methods diverge faster than their square-root counterparts. It is clearly seen from Fig.~2(d) where  the filters' breakdown values are indicated. Thus, new conventional MATLAB-based EKF-5DCKF in Algorithm~1 has the same numerical robustness to roundoff as its previously published NIRK-based variant in~\cite[Alg.~1]{2020:IETSonar:KulikovaKulikov}, but is is much faster; see Fig.~2(c). Secondly, the square-root algorithms are more robust to roundoff than the conventional implementations. More importantly, the numerical stability of the novel MATLAB-based square-root EKF-5DCKF methods and previously published NIRK-based EKF-5DCKF implementations is the same, but the new MATLAB-based methods are faster and much easier to implement and utilize for solving practical problems. Finally, having compared Fig.~1 and Fig.~2, we conclude that the EKF-5DCKF filtering strategy seems to be a more accurate and numerically stable. It is interesting to see that the square-root EKF-5DCKF methods maintain good estimation quality until $\delta = 10^{-12}$ that is in contrast to $\delta = 10^{-11}$ for the EKF-UKF methodology in our test problem.

\section{Concluding remarks} \label{conclusion}

In this paper, a general computational scheme for implementing both the mixed-type  EKF-UKF and EKF-5DCKF estimators are designed with the use of accurate MATLAB ODEs solvers that allow for reducing the discretization error occurred in any {\it continuous-discrete} filtering method. The new filters are simple to implement, they maintain the required accuracy tolerance pre-defined by users at the propagation step with the reduced computational time compared to other estimators of that type published previously. Additionally, their square-root implementations suggested here are robust with respect to roundoff. All this makes the new estimators attractive for solving practical applications.

There are many questions that are still open for a future research. As mentioned in the Introduction, the 5D-CKF moment differential equations have to be derived. This paves a way for designing an accurate 5D-CKF implementation framework with a built-in discretization error control and the reduced moment approximation error compared to the mixed-type EKF-5DCKF suggested in this paper. Besides, it will allow for deriving a simple and elegant 5D-CKF implementation strategy based on the use of MATLAB ODEs solvers. Next, the {\it information}-type filtering methods have their own benefits compared to the {\it covariance}-type algorithms discussed in this work. To the best of the authors' knowledge, the information-type filtering still does not exist for the 5D-CKF estimator. This is an interesting problem for a future research as well. Finally, an alternative approach for designing the square-root filtering methods is based on singular value decomposition (SVD); e.g., see a recent study in~\cite{2020:Automatica:Kulikova}. The spectral methods have the same numerical robustness with respect to roundoff as Cholesky-based square-root algorithms, but they provide users with an extra information about filtering process, e.g., they are capable to oversee the eigenfactors while estimation process, i.e. the singularities may be revealed as soon as they appear~\cite{1988:Oshman,1983:Ham}. Besides, as shown in~\cite{2019:Leth}, eigenvectors used to constructing the sigma points enhance estimation quality of continuous-discrete UKF methods. Unfortunately, it is presently unknown how to derive the SVD-based EKF moment differential equations for implementing the related mixed-type EKF-UKF and EKF-5DCKF estimators. The same problem exists for the SVD-based variant of the UKF and CKF moment differential equations. This also prevents the derivation of accurate and stable MATLAB and SVD-based UKF methods as well as the 3D-CKF and 5D-CKF estimators. Their derivation could be an area for a future research.

%%%%%%%%%%%%%%%%%%%%%%%%%%%%%%%%%%%%%%%%%%%%%%%%%%%%%%%%%%%%%%%%%%%%%%%%5
\section*{Acknowledgments}

The authors acknowledge the financial support of the Portuguese FCT~--- \emph{Funda\c{c}\~ao para a Ci\^encia e a Tecnologia}, through the projects UIDB/04621/2020 and UIDP/04621/2020 of CEMAT/IST-ID, Center for Computational and Stochastic Mathematics, Instituto Superior T\'ecnico, University of Lisbon.
The authors are also thankful to the anonymous referees for their valuable remarks and comments on this paper.

%%%%%%%%%%%%%%%%%%%%%%%%%%%%%%%%%%%%%%%%%%%%%%%%%%%%%%%%%%%%%%%%%%%%%%%%5

\section*{Appendix}

\subsection*{A. Conventional MATLAB-based Time Update filtering schemes}

To present the implementation methodology of the time update steps in the most general way, we denote the utilized solver by \verb"odesolver" in all algorithms presented in this paper and propose a general MATLAB-based time update implementation scheme.
The values $\epsilon_a$ and $\epsilon_r$ are, respectively, absolute and relative tolerances given by users for solving the given ODEs in an accurate and automatic way.  Set the ODEs solver's options with given tolerances $\epsilon_a$, $\epsilon_r$ by \verb"options = odeset(`AbsTol',"$\epsilon_a$,\verb"`RelTol',"$\epsilon_r$).

We start with the discussion of the filters' time update step implementation way.
As can be seen, the second set of the EKF MDEs to be solved is matrix-form equation~\eqref{eq2.2}. One needs to reshape it into a vector-form to be able to use the built-in MATLAB ODEs solvers at the filters' time update step. When the integration is completed, one needs to reshape the resulting solution again in order to recover the state $\hat x_{k|k-1}$ and the filter covariance $P_{k|k-1}$ obtained. In summary, the MATLAB-based time update step of the continuous-discrete EKF with an automatic discretization error control is based on solving the MDEs in~\eqref{eq2.1}, \eqref{eq2.2} and can be summarized in the form of pseudo-code.

\begin{codebox}
\zi \!\!\!\!\!\!\!\!\!\!\textsc{Time Update}: Given $\hat x_{k-1|k-1}$ and $P_{k-1|k-1}$, integrate on $[t_{k-1},t_{k}]$
\li \>Form matrix $XP_{k-1|k-1} = [\hat x_{k-1|k-1}, P_{k-1|k-1}]$;
\li \>Reshape into vector $x^{(0)}_{k-1} = XP_{k-1|k-1}\verb"(:)"$;
\li \>$x_{k|k-1}\leftarrow \verb"odesolver"(\verb"@MDEs",[t_{k-1},t_k],x^{(0)}_{k-1},$
\li \>$\ldots \verb"options",\verb"@drift",\verb"@Jacobian",G,Q,n)$;
\li \>Reshape $XP_{k|k-1} \leftarrow \verb"reshape"(x_{k|k-1}^{(\verb"end")},n,n+1)$ at $t_k$;
\li \>Recover state $\hat x_{k|k-1} = XP_{k|k-1}\verb"(:,1)"$;
\li \>Recover covariance $P_{k|k-1} = XP_{k|k-1}\verb"(:,2:n+1)"$.
\zi \!\!\!\!\!\!\!\!\!\!\textsc{Get the predicted values}: $\hat x_{k|k-1}$ and $P_{k|k-1}$.
\end{codebox}

As can be seen, in the above algorithm, we calculate and work with the entire error covariance matrices $P_{k|k-1}$. Certainly, taking into account a symmetric fashion of any covariance matrix, the size of the MDE system to be solved can be reduced from $n+n^2$ to $n + n(n+1)/2$ equations. However, it requires the full set of such equations to be written done and, then, the symmetric entries of the matrix $P_{k|k-1}$ should be replaced in all $n^2$ equations by hand. In other words, the user needs first to write done the system of MDEs corresponding to the lower (or upper) triangular part of the covariance matrix in use, but the latter could be time-consuming and difficult for some stochastic systems and even experienced practitioners. When this preliminary work is done the resulting differential equations should be reshaped to a vector fashion for application of the built-in MATLAB ODEs solver chosen by the user. We remark that having derived the requested numerical solution one needs the error covariance matrix (or its upper/lower triangular part depending on the approach implemented) and the mean vector to be recovered from this solution. The main advantage of our implementation is that it does all the computations in automatic mode, that is, with no user's effort. The user does not need to do anything except for coding the drift function of stochastic system at hand and to pass it to our algorithm. We consider that our approach is more convenient to practitioners because it works in the same way for any stochastic system which may potentially arise in practice. From our point of view, the convenience of utilization of any method employed in applied science and engineering is preferable than a negligible reduction in the state estimation time. However, in a particular case of computational device's limited power and memory, the above-mentioned MDE reduction can be implemented by the thoughtful readers, if necessary.

The function for computing the right-hand side of system~\eqref{eq2.1}, \eqref{eq2.2} should be sent to the chosen built-in MATLAB ODEs solver.

\begin{codebox}
\Procname{$[xnew] \leftarrow \proc{MDEs}(t,x,drift,Jac,Q,G,n)$}
\li \>\verb"A = reshape(x,n,n+1); % reshape into matrix"
\li \>\verb"X = A(:,1);           % read-off state"
\li \>\verb"P = A(:,2:n+1);       % read-off covariance"
\li \>\verb"J = feval(Jac,t,X);   % compute Jacobian"
\li \>\verb"X = feval(drift,t,X); % implement eq. (3)"
\li \>\verb"P = J*P + P*J'+G*Q*G' % implement eq. (4)"
\li \>\verb"A = [X, P];           % collect matrix"
\li \>\verb"xnew = A(:);          % reshape into vector"
\end{codebox}

\subsection*{B. Square-root MATLAB-based Time Update filtering schemes}

Within the square-root methods, the Cholesky decomposition is applied to the initial filter error covariance $\Pi_0$. Next, the square-root factors $P^{1/2}_{k|k-1}$ are propagated instead of the entire error covariance matrices $P_{k|k-1}$ in each iteration step. Thus, we obtain the following square-root version of the MATLAB-based time update step.

\begin{codebox}
\zi \!\!\!\!\!\!\textsc{Time Update}: \Comment{\small\textsc{integrate on $[t_{k-1},t_{k}]$ in line with given $\epsilon_a$, $\epsilon_r$}}
\li \>Form $XP_{k-1|k-1} = [\hat x_{k-1|k-1}, P^{1/2}_{k-1|k-1}]$;
\li \>Reshape into vector $x^{(0)}_{k-1} = XP_{k-1|k-1}\verb"(:)"$;
\li \>$x_{k|k-1}\leftarrow \verb"odesolver"(\verb"@MDESR",[t_{k-1},t_k],x^{(0)}_{k-1},\verb"options")$;
\li \>Reshape $XP_{k|k-1} \leftarrow \verb"reshape"(x_{k|k-1}^{(\verb"end")},n,n+1)$ at $t_k$;
\li \>Recover state $\hat x_{k|k-1} = XP_{k|k-1}\verb"(:,1)"$;
\li \>Recover factor $P_{k|k-1}^{1/2} = XP_{k|k-1}\verb"(:,2:n+1)"$;
\zi \!\!\!\!\!\!\!\!\!\!\textsc{Get the predicted values}: $\hat x_{k|k-1}$ and $P^{1/2}_{k|k-1}$.
\end{codebox}

Additionally, the pseudo-code below performs the right-hand side function of the MDE system~\eqref{eq2.1}, \eqref{eq3.13} in its square-root form.

\begin{codebox}
\Procname{$[xnew] \leftarrow \proc{MDESR}(t,x,drift,Jac,Q,G,n)$}
\li \>\verb"L = reshape(x,n,n+1);% reshape into matrix"
\li \>\verb"X = L(:,1);          % read-off state"
\li \>\verb"S = L(:,2:n+1);      % read-off SR factor"
\li \>\verb"J = feval(Jac,t,X);  % compute Jacobian"
\li \>\verb"X = feval(drift,t,X);% implement eq.(3)"
\li \>\verb"A = S\(J*S);         % find  for eq.(19)"
\li \>\verb"M = A + A' + S\(G*Q*G')/(S');"
\li \>\verb"Phi = tril(M,-1) + diag(diag(M)/2);"
\li \>\verb"L = [X, S*Phi];      % implement eq.(19)"
\li \>\verb"xnew = L(:);         % reshape in vector"
\end{codebox}

%\bibliographystyle{model1b-num-names}
%\bibliography{BibTex_Library/books,%
%              BibTex_Library/KFDistribution_Robust,%
%              BibTex_Library/KF_Chandrasekhar,%
%              BibTex_Library/KFMCC_Riccati,%
%              BibTex_Library/KFMCC_Chandrasekhar,%
%              BibTex_Library/KFMCC_Applications,%
%              BibTex_Library/KFDiff_Chandrasekhar,%
%              BibTex_Library/KF_Riccati,%
%              BibTex_Library/CKF_Riccati,%
%              BibTex_Library/UKF_Riccati,%
%              BibTex_Library/EKF_all,%
%              BibTex_Library/KF_Applications,%
%              BibTex_Library/list_Kulikov,%
%              BibTex_Library/list_KF_nonGaussian,%
%              BibTex_Library/Lin_Algebra}

\end{document}